\documentclass[reqno,centertags]{amsart}
\usepackage{amsmath,amsthm,amscd,amssymb}
\usepackage{latexsym}
\newcommand{\bbC}{{\mathbb{C}}}
\newcommand{\bbD}{{\mathbb{D}}}

\newcommand{\bbR}{{\mathbb{R}}}

\newcommand{\bbZ}{{\mathbb{Z}}}

\newcommand{\calB}{{\mathcal{B}}}
\newcommand{\calC}{{\mathcal{C}}}

\newcommand{\calS}{{\mathcal S}}

\newcommand{\dott}{\,\cdot\,}

\newcommand{\lb}{\label}
\newcommand{\f}{\frac}

\newcommand{\ol}{\overline}
\newcommand{\ti}{\tilde  }

\newcommand{\tr}{\text{\rm{Tr}}}

\newcommand{\bi}{\bibitem}

\newcommand{\beq}{\begin{equation}}
\newcommand{\eeq}{\end{equation}}
\newcommand{\ba}{\begin{align}}
\newcommand{\ea}{\end{align}}

\newcounter{smalllist}
\newenvironment{SL}{\begin{list}{{\rm\roman{smalllist})}}{%
\setlength{\topsep}{0mm}\setlength{\parsep}{0mm}\setlength{\itemsep}{0mm}%
\setlength{\labelwidth}{2em}\setlength{\leftmargin}{2em}\usecounter{smalllist}%
}}{\end{list}}

\DeclareMathOperator{\Real}{Re}
\DeclareMathOperator{\Ima}{Im}

\allowdisplaybreaks
\numberwithin{equation}{section}

\newtheorem{theorem}{Theorem}[section]

\newtheorem*{p2.1}{Proposition 2.1}
\newtheorem{proposition}[theorem]{Proposition}
\newtheorem{lemma}[theorem]{Lemma}

\theoremstyle{definition}
\newtheorem{example}[theorem]{Example}

\theoremstyle{remark}
\newtheorem*{remark}{Remark}
\newtheorem*{remarks}{Remarks}

\newcommand{\abs}[1]{\lvert#1\rvert}

\begin{document}
\title{Poisson Brackets of Orthogonal Polynomials}
\author[M.~J.~Cantero and B.~Simon]{Mar\'ia Jos\'e Cantero$^1$ and Barry Simon$^2$}

\thanks{$^1$ Department of Applied Mathematics, University of Zaragoza, 50018 Zaragoza,
Spain. E-mail: mjcante@unizar.es. Supported in part by a research
grant from the Ministry of Education and Science of Spain, MTM2005-08648-C02-01,
and by the Arag\'on government, Project E-64 of DGA}

\thanks{$^2$ Mathematics 253-37, California Institute of Technology, Pasadena, CA 91125.
E-mail: bsimon@caltech.edu. Supported in part by NSF Grant DMS-0140592 and
U.S.--Israel Binational Science Foundation (BSF) Grant No.\ 2002068}

\date{March 9, 2007}

\begin{abstract} For the standard symplectic forms on Jacobi and
CMV matrices, we compute Poisson brackets of OPRL and OPUC, and relate
these to other basic Poisson brackets and to Jacobians of basic changes of
variable.
\end{abstract}

\maketitle

\section{Introduction} \lb{s1}

It has been known since the discoveries of Flaschka \cite{Fla1} and Moser \cite{Mo75}
concerning the finite Toda lattice as a completely integrable system that finite
Jacobi matrices of fixed trace support a natural symplectic form. Explicitly, if
\begin{equation} \lb{1.1}
J= \begin{pmatrix}
b_1 & a_1 & {}  & {} \\
a_1 & b_2 &  a_2  \\
{} & \ddots  & \ddots  \\
{} &  {} & a_{N-1} &  b_N
\end{pmatrix}
\end{equation}
on the manifold where $a_j >0$ and 
\begin{equation} \lb{1.2}
\sum_{j=1}^N b_j = c
\end{equation}
the nonzero Poisson brackets of the $a$'s and $b$'s are
\begin{alignat}{2}
\{b_k, a_k\} &= -\tfrac14\, a_k \qquad && k=1, \dots, N-1 \lb{1.3} \\
\{b_k, a_{k-1}\} &= \tfrac14\, a_{k-1} \qquad && k=2, \dots, N \lb{1.4}
\end{alignat}
For the basic definitions of symplectic manifolds and Poisson brackets,
see Deift \cite{Deift96}.

The nonzero pairings are neighboring elements in $J$. As one would guess,
there are differing conventions, so that sometimes the $\f14$ in \eqref{1.3}
and \eqref{1.4} are $\f12$ and sometimes one takes $\{\,\, , \,\,\}$ to be
the negative of our choice.

It is easy to check that $\sum_{j=1}^N b_j$ commutes with all $a$'s and $b$'s
as is necessary for this to define a bracket on the manifold where \eqref{1.2}
holds. It is also easy to see that the form defined by \eqref{1.3}/\eqref{1.4} is
nondegenerate and closed so there is an underlying symplectic form; see
the end of Section~\ref{s3}.

If one takes
\begin{equation} \lb{1.5}
H=2 \, \tr (J^2) =2\sum_{j=1}^N b_j^2 + 4\sum_{j=1}^{N-1} a_j^2
\end{equation}
then the Hamiltonian flow,
\begin{align}
\Dot{b}_k &= \{H,b_k\} = 2(a_k^2 -a_{k-1}^2) \lb{1.6} \\
\Dot{a}_k &= \{H,a_k\} = a_k (b_{k+1} -b_k) \lb{1.7}
\end{align}
is the Toda flow in Flaschka form \cite{Fla1}. It is completely integrable;
indeed, (see \eqref{1.11} below)
\begin{equation} \lb{1.8}
\{\tr(J^n), \tr (J^m)\}=0
\end{equation}
for all $n,m$.

If $d\rho$ is the spectral measure for $J$ and $(1,0,\dots, 0)^t$, then
\begin{equation} \lb{1.9}
d\rho = \sum_{j=1}^N \rho_j \delta_{x_j}
\end{equation}
The $\{\rho_j\}_{j=1}^N$ and $\{x_j\}_{k=1}^N$ are not independent; indeed,
by \eqref{1.2},
\begin{equation} \lb{1.10}
\sum_{j=1}^N x_j =c \qquad
\sum_{j=1}^N \rho_j =1
\end{equation}

Fundamental to developments are the following Poisson brackets:
\begin{equation} \lb{1.11}
\{x_j, x_k\} =0\qquad
\{x_j,\rho_k\} = \tfrac12\, [\delta_{jk} \rho_j - \rho_j\rho_k]
\end{equation}

The analog of $\{x_j, x_k\}=0$ for periodic Toda chains is due to Flaschka
\cite{Fla1} and \eqref{1.11} appeared (implicitly) first in Moser \cite{Mo75}.

Our original motivation was to  understand the analog of this for CMV matrices.
Along the way, we realized that Poisson brackets of orthogonal polynomials were
not previously studied, even in the Jacobi case, so we decided to discuss both
cases.

We also felt that the centrality of \eqref{1.11}, which is partly known to experts,
was not always so clear, so we also decided to say something about that. Indeed,
one of our initial goals was to prove the analog of $\{x_i, x_j\}=0$
for the periodic case of OPUC whose previous proofs \cite[Section~11.11]{OPUC2}
and \cite{NenDiss,Nen05} were involved. Hence one application of \eqref{1.11} we
discuss is going from \eqref{1.11} to the periodic case.

Since $\sum_{j=1}^N x_j$ and $\sum_{j=1}^N \rho_j$ are constant, $\{x_j, \rho_k\}$
must sum to zero over $j$ or $k$, as can be checked in \eqref{1.11} since $\sum_{j=1}^N
\rho_j =1$.

We note one curiosity of \eqref{1.11}, namely, $\{x_j,\rho_k\}$ is symmetric in
$j$ and $k$. We will eventually also find
\begin{equation} \lb{1.11a}
\{\rho_j,\rho_k\} = \f{\rho_j\rho_k}{x_j-x_k} - \sum_{m\neq j}
\f{\rho_j\rho_k\rho_m}{x_j-x_m} + \sum_{m\neq k} \f{\rho_j\rho_k\rho_m}{x_k-x_m}
\end{equation}
but will make no use of this complicated formula here. It is needed if one wants to
find angle variables for the Toda flows (see \cite{KN2,GNprep}).

Simultaneous to our work, Gekhtman--Nenciu \cite{GNprep} were studying Poisson brackets
of Carath\'eodory functions, which we will see is closely related to our work. Indeed,
a key to our progress was learning of this work and the earlier paper of
Faybusovich--Gekhtman \cite{FG2000}. We also mention not unrelated earlier papers 
of Kako--Mugibayashi \cite{KM78,KM79} and Kulish \cite{Kul}. 

For the monic orthogonal polynomials on the real line (OPRL), $P_n(x)$, our basic
result is that for $n=1,2,\dots, N$,
\begin{align}
\{P_n(x), P_n(y)\} &= \{P_{n-1}(x), P_{n-1}(y)\} = 0 \lb{1.12} \\
2\{P_n(x), P_{n-1}(y)\} &= \biggl[ \f{P_n(x) P_{n-1}(y) - P_n(y) P_{n-1}(x)}{x-y}\biggr]
- P_{n-1}(x) P_{n-1}(y) \lb{1.13}
\end{align}

Notice the occurrence of the Bezoutian (see \cite{HF89}) of $P_n$ and $P_{n-1}$ in the first
term on the right of \eqref{1.13}. \eqref{1.12} and \eqref{1.13} will follow by a direct
induction in $n$.

There is a symmetry between $P_n, P_{n-1}$ and $P_n,Q_n$ (where $Q_n$ is the second kind
polynomial (of degree $n-1$) for $P_n$), for $P_N$ is a determinant of $z-J$ and $P_{N-1}$
(resp.\ $Q_N$) are the determinants obtained by removing the last and rightmost row
and column (resp.\ first and leftmost row and column). This will immediately lead to
\begin{align}
\{P_n (x), P_n(y)\} &= \{Q_n(x), Q_n (y)\} =0 \lb{1.14} \\
2\{P_n(x), Q_n(y)\} &= - \biggl[ \f{P_n(x) Q_n(y) - P_n(y) Q_n(x)}{x-y}\biggr]
+ Q_n(x) Q_n(y) \lb{1.15}
\end{align}
The sign changes because inverting order changes the sign of $\{\,\, , \,\,\}$.
All these calculations appear in Section~\ref{s2}. In \cite{FG2000}, several
different Poisson brackets on polynomials are considered, including one that
is essentially equivalent to \eqref{1.15}.

In Section~\ref{s3}, we derive $\{x_j,x_k\}=0$ from $\{P_n(x), P_n(y)\}=0$,
$\{x_j, \rho_k\}$ from $\{P_n(x), Q_n(y)\}$ and $\{\rho_j, \rho_k\}$
from $\{Q_n(x), Q_n(y)\} =0$.

Looked at carefully, if we argued directly for $P,Q$ by induction, we, in essence,
use coefficient stripping to relate $P_N,Q_N$ for $\{a_j\}_{j=1}^{N-1}\cup
\{b_j\}_{j=1}^N$ to $P_{N-1}, Q_{N-1}$ for $\{a_j\}_{j=2}^{N-1}\cup
\{b_j\}_{j=2}^{N-1}$. The analog for orthogonal polynomials on the unit
circle (OPUC) is a new coefficient stripping relation of first and
second kind paraorthogonal polynomials, which is new. We present this
in Section~\ref{s4}, denoting these first and second kind OPUC by
$P_N$ and $Q_N$.

In Section~\ref{s5}, we use the natural symplectic form for OPUC: in terms of
Verblunsky coefficients,
\begin{equation} \lb{1.16}
\{\alpha_j, \alpha_k\} =0 \qquad
\{\alpha_j, \bar\alpha_k\} = -i\rho_j^2 \delta_{jk}
\end{equation}
This was introduced by Nenciu--Simon \cite{NS}. We will find
\begin{align}
\{P_n(z), P_n(w)\} &= 0 = \{Q_n(z), Q_n(w)\} \lb{1.17} \\
\{P_n(z), Q_n(w)\} &=  -\f{i}{2} \biggl[ (P_n(z) Q_n(w) - P_n(w) Q_n(z))
\biggl( \f{z+w}{z-w}\biggr) \notag \\
&\qquad \qquad \qquad
- Q_n(z) Q_n(w) + P_n (z) P_n(w)\biggr] \lb{1.18}
\end{align}
using induction and the coefficient stripping of Section~\ref{s4}. $P_n, Q_n$
will depend on a parameter $\beta$. $\{\alpha_j\}_{j=0}^{n-2} \to \{-\alpha_j
\}_{j=0}^{n-2}$ which preserves OPUC Poisson brackets and $\beta \to -\beta$
interchanges $P_n$ and $Q_n$, so the right side of \eqref{1.19} has to
be antisymmetric under interchange of $P_n$ and $Q_n$, as it is.

In Section~\ref{s6}, we derive the analogs of \eqref{1.11}. In \eqref{1.11},
the $x$'s and $\rho$'s are global variables fixed by $x_1 < x_2 < \dots$.
For OPUC, the eigenvalues are only locally well defined, namely,
\begin{equation} \lb{1.19}
d\mu = \sum_{j=1}^N \mu_j \delta_{z_j}
\end{equation}
where
\begin{equation} \lb{1.20}
z_j = e^{i\theta_j}
\end{equation}
The analogs of \eqref{1.10} are
\begin{equation} \lb{1.21}
\sum_{j=1}^N \theta_j = c \qquad
\sum_{j=1}^N \mu_j =1
\end{equation}
and of \eqref{1.11},
\begin{equation} \lb{1.22}
\{\theta_j, \theta_k\}=0 \qquad
\{\theta_j, \mu_k\} = \mu_j \delta_{jk} - \mu_j \mu_k
\end{equation}

The analog of $\{\theta_j, \theta_k\}=0$ for the periodic case appeared first in
\cite{NS}, and $\{\theta_j,\theta_k\}=0$ is from \cite{NenDiss,Nen05}.
The formula for $\{\theta_j, \mu_k\}$ is due to Killip--Nenciu \cite{KN2}.

In Sections~\ref{s7}--\ref{s9}, we discuss three applications of the fundamental
relations \eqref{1.11} and \eqref{1.22}: to computations of Jacobians of the
maps $(a,b)\to (x,\rho)$ and $(\alpha)\to (\theta,\mu)$, to the exact solution
of the flows generated by $G(x)$ or $G(\theta)$, and to the periodic case.
In Sections~\ref{s9new} and \ref{s10new}, we discuss the differential equations
induced on the monic OPRL and OPUC. Sections~\ref{s10} and \ref{s11} compute
more Poisson brackets for OPRL and OPUC, respectively.

\medskip
It is a pleasure to thank Rafael Hern\'andez Heredero and Irina Nenciu for useful discussions.
M.~J.~Cantero would like to thank Tom Tombrello and Gary Lorden for the
hospitality of Caltech.

\section{Poisson Brackets of First and Second Kind OPRL} \lb{s2}

Our goal in this section is to prove \eqref{1.12} through \eqref{1.15}. We note that since 
decreasing $N$ to $N-1$ does not change the Poisson brackets for $\{a_k\}_{k=1}^{N-2}, 
\{b_j\}_{j=1}^{N-1}$ nor the $\{P_j\}_{j=1}^{N-1}$, we can suppose $n=N$ in proving 
\eqref{1.12}--\eqref{1.15} but then increase $N$ in the induction step. 

\begin{proposition} \lb{P2.1} \eqref{1.12} and \eqref{1.13} holds for $n=1$.
\end{proposition}

\begin{proof} For $N=1$, $b_1=c$, $P_0(x)=1$, $P_1(x)=x-c$, and all Poisson brackets are
zero. Thus, we need only show the RHS of \eqref{1.13} is zero, that is,
\[
\f{(x-c)-(y-c)}{x-y} -1 =0
\]
which is obvious.
\end{proof}

\begin{lemma}\lb{L2.2}
\begin{equation} \lb{2.1}
\{a_n^2, P_n(w)\} = -\tfrac12\, a_n^2 P_{n-1}(w)
\end{equation}
\end{lemma}

\begin{proof} Since
\begin{equation} \lb{2.2}
P_{j+1}(x) = (x-b_{j+1}) P_j(x) - a_j^2 P_{j-1}(x)
\end{equation}
we have inductively that $P_j$ is a function of $\{b_\ell\}_{\ell=1}^j \cup
\{a_\ell\}_{\ell=1}^{j-1}$ as also follows from \eqref{2.8} below. Thus,
we have $\{a_n^2, P_{n-1}(w)\}=\{a_n^2, P_{n-2}(w)\} =0$ so, using
\eqref{2.2} for $j+1=n$,
\begin{align*}
\{a_n^2, P_n(w)\} &= \{a_n^2, -b_n\} P_{n-1}(w) \\
&= -\tfrac12\, a_n^2 P_{n-1}(w)
\qedhere
\end{align*}
\end{proof}

\begin{theorem}\lb{T2.3} \eqref{1.12} and \eqref{1.13} hold for all $n$.
\end{theorem}

\begin{proof} Proposition~\ref{P2.1} is the result for $n=1$. So, by induction,
we can suppose \eqref{1.12} and \eqref{1.13} for $j\leq n$ and need only prove it
for $n+1$. By \eqref{2.2} for $j=n$ and the facts that $\{P_{n-1}(x), P_{n-1}(y)\}
= \{P_n(x), P_n(y)\}=0$ and $P_n$ independent of $a_n$, $a_{n+1}$ and $P_{n-1}$
independent of $b_n$ and $b_{n+1}$, only the cross terms enter and
\[
- \{P_{n+1}(x), P_{n+1}(y)\} = \{(x-b_{n+1}) P_n(x), a_n^2 P_{n-1}(y)\} - (x\leftrightarrow y)
\]
where $(x\leftrightarrow y)$ is like the first term but $x$ and $y$ are reversed.

This is a sum of three terms: $t_1, t_2, t_3$ where
\begin{align}
2t_1 &= 2(x-b_{n+1}) a_n^2 \{P_n(x), P_{n-1}(y)\} - (x\leftrightarrow y) \notag \\
&= (x-y) a_n^2 \biggl[ \f{(P_n(x) - P_{n-1}(y) - P_n (y) P_{n-1}(x))}{x-y}
- P_{n-1}(x) P_{n-1}(y)\biggr] \lb{2.3}
\end{align}
by the induction hypothesis and the symmetry of $\{P_n(y), P_{n-1}(y)\}$ under
$x\leftrightarrow y$. Next,
\begin{align}
2t_2 &= 2\{xP_n(x), a_n^2\} P_{n-1}(y) - (x\leftrightarrow y) \notag \\
&= a_n^2 (x-y) P_{n-1}(x) P_{n-1}(y) \lb{2.4}
\end{align}
by Lemma~\ref{L2.2}. Finally,
\begin{align}
2t_3 &= -2\{b_{n+1}, a_n^2\} (P_n(x) P_{n-1}(y) - (x-y)) \notag \\
&= -a_n^2 (P_n(x) P_{n-1}(y) - P_n(y) P_{n-1}(x)) \lb{2.5}
\end{align}
There is a fourth term involving $\{b_{n+1}, P_{n-1}(y)\}$, but this is zero.

Clearly, $t_1 + t_2 + t_3 =0$ since \eqref{2.5} cancels the first term in \eqref{2.3}, and
\eqref{2.4} the second term. This proves that $\{P_{n+1}(x), P_{n+1}(y)\}=0$.

Since $\{P_n(x), P_n(y)\}=0$ and $\{b_{n+1}, P_n(y)\}=0$ (since $P_n(y)$ only depends on
$\{a_j\}_{j=0}^{n-1}$), \eqref{2.2} implies
\begin{align}
\{P_{n+1}(x), P_n(y)\} &= \{-a_n^2 P_{n-1}(x), P_n(y)\} \notag \\
&= -a_n^2 \{P_{n-1}(x), P_n(y)\} - P_{n-1}(x) \{a_n^2, P_n(y)\} \lb{2.6}
\end{align}
The first term is evaluated by induction and the second by Lemma~\ref{L2.2}. The Lemma~\ref{L2.2}
term is $+\f12 a_n^2 P_{n-1}(x) P_{n-1}(y)$ which cancels one term in $\{P_{n-1}(x),
P_n(y)\}$. Thus
\begin{align}
2\{P_{n+1}(x), P_n(y)\}
&= a_n^2 \biggl[ \f{P_n(y) P_{n-1}(x) - P_n(x) P_{n-1}(y)}{y-x}\biggr] \notag \\
&= \biggl[ -\f{P_n(y) (P_{n+1}(x) - (x-b_{n+1}) P_n(x)) - (x\leftrightarrow y)}{y-x}\biggr] \lb{2.7} \\
&= \f{(P_{n+1}(y) P_n(y) - P_{n+1}(y) P_n(x))}{x-y} - P_n(x) P_n(y) \notag
\end{align}
proving the required formula for $n+1$. In the above, \eqref{2.7} comes from \eqref{2.2}.
\end{proof}

The monic second kind polynomials, $Q_n$, can be defined as follows: We note that if
$J(b_1, \dots, b_n; a_1, \dots, a_{n-1})$ is the matrix \eqref{1.1} with $N=n$, then
\begin{equation} \lb{2.8}
P_n(x) = \det (x-J(b_1, \dots, b_n; a_1, \dots, a_{n-1}))
\end{equation}
as is well known. $Q_n$ is then defined by removing the top row and left column, that is,
\begin{equation} \lb{2.9}
Q_n(x) = \det (x-J(b_2, \dots, b_n; a_2, \dots, a_{n-1}))
\end{equation}

Notice that $J(b_1, \dots, b_n; a_1, \dots, a_{n-1})$ and $J(b_n, b_{n-1}, \dots,
b_1; a_{n-1}, \dots, a_1)$ are unitarily equivalent under the unitary $U\delta_j =
\delta_{n+1 -j}$, $j=1, \dots, n$. This shows making the $a,b$ dependence explicit
(we include redundant variables; $Q_n$ is independent of $b_1$ and $a_1$):
\begin{align*}
P_n(x; b_1, \dots, b_n; a_1, \dots, a_{n-1}) &= P_n (x; b_n, \dots, b_1; a_{n-1}, \dots, a_1) \\
Q_n(x; b_1, \dots, b_n; a_1, \dots, a_{n-1}) &= P_{n-1} (x; b_n, \dots b_1; a_{n-1}, \dots, a_1)
\end{align*}
Taking into account that reversing order changes the signs in \eqref{1.13} and \eqref{1.15}, we
see that $\{\,\, , \,\,\}$ changes sign, that is, \eqref{1.13} implies \eqref{1.15}, and we
have

\begin{theorem}\lb{T2.5} \eqref{1.14} and \eqref{1.15} hold for all $n$.
\end{theorem}

Instead of using induction from $P_n, P_{n-1}$ to $P_{n+1}, P_n$ and the above reversal
argument, we can directly use induction from $P_n,Q_n$ to $P_{n+1}, Q_{n+1}$ using
\begin{equation} \lb{2.10}
Q_{n+1}(x; b_1, \dots, b_{n+1}; a_1, \dots, a_n) = P_n (x; b_2, \dots b_{n+1}; a_2, \dots, a_n)
\end{equation}
and
\begin{equation} \lb{2.11}
\begin{split}
&P_{n+1}(x; b_1, \dots, b_{n+1}; a_1, \dots, a_n) \\
 &\, \, = (x-b_1) P_n (x; b_2, \dots, b_{n+1}; a_2, \dots, a_n)
- a_1^2 Q_n (x; b_2, \dots, b_{n+1}; a_2, \dots, a_n)
\end{split}
\end{equation}
that is, coefficient stripping.

\section{Fundamental Poisson Brackets for OPRL} \lb{s3}

Our goal in this section is to prove \eqref{1.11a} and, most importantly,
\eqref{1.11}. In the calculations below, we will compute Poisson brackets
of functions of $a,b$ and free parameter(s), $x$ (e.g., $P_n(x)$) and
then set $x$ to a value that is a function, $h$, of $a,b$. The order of
operations is important, that is, $\{\cdot, \cdot\}$ first, only then
evaluation. We will denote this by $\left. \{f,g\}\right|_{x=h}$.

\begin{proposition}\lb{P3.1}
\begin{equation} \lb{3.1}
\left. \{P_n(x), P_n(y)\}\right|_{x=x_j,\, y=x_k} =
\{x_j, x_k\} \prod_{\ell\neq j} (x_j-x_\ell)
\prod_{\ell\neq k} (x_k-x_\ell)
\end{equation}
In particular,
\begin{equation} \lb{3.2}
\{P_n(x), P_n(y)\}=0\Leftrightarrow \{x_j, x_k\}=0\qquad\text{for all } j,k
\end{equation}
\end{proposition}

\begin{proof} Since $P$ is monic and its zeros are simple and \eqref{2.8} holds,
\begin{equation} \lb{3.3}
P_n(x) = \prod_j (x-x_j)
\end{equation}
Thus, by Leibnitz's rule,
\begin{equation} \lb{3.4}
\{P_n(x), P_n(y)\} =\sum_{p,q} \{x_p, x_q\} \prod_{\ell\neq p} (x-x_\ell)
\prod_{\ell\neq q} (y-x_\ell)
\end{equation}
Setting $x=x_j$, $y=x_k$, all terms with $p\neq j$ or $q\neq k$ vanish and we
obtain \eqref{3.1}. \eqref{3.1} trivially implies $\Rightarrow$ in \eqref{3.2} 
and \eqref{3.4} implies $\Leftarrow$ in \eqref{3.2}.
\end{proof}

\begin{proposition}\lb{P3.2} We have that
\begin{equation} \lb{3.5}
Q_n(x) =\sum_{j=1}^n \rho_j \prod_{k\neq j} (x-x_k)
\end{equation}
\end{proposition}

\begin{proof} The $m$-function is defined by
\begin{equation} \lb{3.6}
m_n(z) =\sum_{j=1}^n \f{\rho_j}{x_j-z}
\end{equation}
which is just
\begin{equation} \lb{3.7}
m_n(z)=\langle\delta_1, (J(b_1, \dots, b_n; a_1, \dots, a_{n-1})-z)^{-1}\delta_1\rangle
\end{equation}
which, by Cramer's rule and \eqref{2.8}/\eqref{2.9}, implies
\begin{equation} \lb{3.8}
m_n(z) = -\f{Q_n(z)}{P_n(z)}
\end{equation}
Thus, by \eqref{3.3},
\begin{align*}
Q_n(z) &= - m_n(z) P_n(z) \\
&= \sum_{j=1}^n \f{\rho_j (\prod_k (z-x_k))}{(z-x_j)}
\end{align*}
which is \eqref{3.5}.
\end{proof}

\begin{proposition} \lb{P3.3}
\begin{equation} \lb{3.9}
\left. \{P_n(x), Q_n(y)\}\right|_{x=x_j,\, y=x_k} = -\{x_j, \rho_k\}
\prod_{\ell\neq j} (x_j-x_\ell) \prod_{\ell\neq k} (x_k -x_\ell)
\end{equation}
This implies that
\begin{equation} \lb{3.10}
\eqref{1.15} \Leftrightarrow \{x_j, \rho_k\} = \tfrac12\, [\delta_{jk} \rho_j - \rho_j \rho_k]
\end{equation}
\end{proposition}

\begin{proof} Since we have proven that $\{x_j, x_k\}=0$, \eqref{3.3} and \eqref{3.5} imply
that
\begin{equation} \lb{3.11}
\{P_n(x), Q_n(y)\} =-\sum_{p,q} \, \{x_p, \rho_q\} \prod_{\ell\neq p} (x-x_\ell)
\prod_{\ell\neq q} (y-x_\ell)
\end{equation}
which immediately implies \eqref{3.9}.

In \eqref{1.15}, we have, by \eqref{3.5}, that
\begin{equation} \lb{3.12}
\left. Q_n(x) Q_n(y)\right|_{x=x_j, \, y=x_k} = \{\rho_j \rho_k\} \prod_{\ell\neq j}
(x_j -x_\ell) \prod_{\ell\neq k} (x_k-x_\ell)
\end{equation}
The Bezoutian term in \eqref{1.15} is $0$ if $j\neq k$, but is slightly subtle
if $j=k$ because then the formal expression is $0/0$. We thus take limits since,
of course, $\{P_n(x), Q_n(y)\}$, as a polynomial in $x$ and $y$, is continuous.
Since $P_n(x_k)=0$, for $x\neq x_k$,
\begin{equation} \lb{3.13}
\f{P_n(x) Q_n(x_k) - P_n(x_k) Q_n(x)}{x-x_k} = \biggl[ \, \prod_{\ell\neq k}
(x-x_\ell)\biggr] \rho_k \biggl[\, \prod_{\ell\neq k} (x-x_\ell)\biggr]
\end{equation}
Thus
\begin{equation} \lb{3.14}
\lim_{\substack{ x\to x_j \\ x\neq x_k}}\, \text{LHS of \eqref{3.13}} =
\rho_k \biggl[\, \prod_{\ell\neq k} (x_k -x_\ell)\biggr]^2 \delta_{jk}
\end{equation}

Thus \eqref{1.15} implies the right side of \eqref{3.10}. By \eqref{3.11}, the 
right side of \eqref{3.10} implies \eqref{1.15}. 
\end{proof}

\begin{remark} The symmetry of $\{x_j, \rho_k\}$ under $j\leftrightarrow k$
is equivalent to the symmetry of the right side of \eqref{1.15} under
$x\leftrightarrow y$.
\end{remark}

\begin{proposition} \lb{P3.4} Let $j\neq k$. Then
\begin{align}
\{\rho_j, \rho_k\} &= \sum_{m\neq j} \f{\rho_j \{x_m, \rho_k\} + \rho_m \{x_j,\rho_k\}}{x_j-x_m}
-\sum_{m\neq k} \f{\rho_k \{x_m,\rho_j\} + \rho_m \{x_k,\rho_j\}}{x_k-x_m} \lb{3.15} \\
&= \f{\rho_j \rho_k}{x_j-x_k} -\sum_{m\neq j} \f{\rho_j \rho_k\rho_m}{x_j-x_m}
+\sum_{m\neq k} \f{\rho_j\rho_k\rho_m}{x_k-x_m} \lb{3.16}
\end{align}
\end{proposition}

\begin{proof} \eqref{3.15} comes from $\{Q_n(x), Q_n(y)\}=0$. Explicitly,
\[
\f{\{Q_n(x), Q_n(y)\}_{x=x_j,\, y=x_k}}
{\prod_{\ell \neq j} (x_j-x_\ell) \prod_{\ell\neq k} (x_k-x_\ell)} =
\{\rho_j,\rho_k\} - \text{RHS of \eqref{3.15}}
\]
by a straightforward but tedious calculation. \eqref{3.16} then follows by using
\eqref{1.11} for $\{x_p,\rho_q\}$.
\end{proof}

\begin{remark} There is a sense in which \eqref{1.11} implies \eqref{1.11a} since
\eqref{1.11} implies $\{P_n(x), P_n(y)\}=0$ which means $\{Q_n(x), Q_n(y)\}=0$
(since $Q_n$ is a $P_{n-1}$ for suitable Jacobi parameters), and this yields
\eqref{3.15}.
\end{remark}

Finally, we want to translate the basic Poisson brackets relations \eqref{1.13}/\eqref{1.14}
and \eqref{1.11} into assertions about the two-form which defines the symplectic structure
underlying the Poisson brackets. Recall a symplectic structure is defined by a two-form, $\omega$,
that is, an antisymmetric functional on tangent vectors with the requirement that $\omega$
is closed (which is equivalent to the Jacobi identity for the Poisson bracket; see, e.g.,
Deift \cite{Deift96}) and which is nondegenerate; that is, for all tangent vectors
$v\neq 0$ at $p\in M$, there is $\ti v$ at $p$ so $\omega(v,\ti v)\neq 0$.
Nondegeneracy implies $d=\dim(M)$ is even.

Given such a two-form and function, $f$, the Hamiltonian vector field $H_f$
is defined by
\begin{equation} \lb{3.17}
\omega(H_f, v) = \langle df, v\rangle
\end{equation}
for all tangent vectors $v$. Here $df$ is the one-form $\sum_{i=1}^d
\f{\partial f}{\partial x_i}dx_i$ as usual, and $\langle \cdot, \cdot\rangle$
is the linear algebra pairing of vectors and their duals.

The Poisson bracket associated to $\omega$ is then
\begin{equation} \lb{3.18}
\{f,g\}= H_f (g)
\end{equation}

\begin{proposition} \lb{P3.5} Let $M$ be a symplectic manifold of dimension $d=2\ell$ and
suppose there is a local coordinate system $\{x_j\}_{j=1}^\ell$, $\{y_j\}_{j=1}^\ell$, so
that
\begin{equation} \lb{3.19}
\omega = \sum_{i,j=1}^\ell W_{ij} \, dx_i \wedge dy_j + \sum_{i,j=1}^\ell
U_{ij}\, dx_i \wedge dx_j
\end{equation}
{\rm{(}}i.e., no $dy_i \wedge dy_j$ terms{\rm{)}} with $U$\! antisymmetric. Then
\begin{SL}
\item[{\rm{(i)}}] Nondegeneracy of $\omega$ is equivalent to invertibility of the
$\ell\times\ell$ matrix $(W_{ij})_{1\leq i,j\leq \ell}$. We will denote its inverse
by $(W^{-1})_{ij}$.

\item[{\rm{(ii)}}]
\begin{equation} \lb{3.20}
\omega^\ell= \ell!\,\det(W) dx_1 \wedge dy_1 \wedge dx_2 \wedge \dots \wedge dy_\ell
\end{equation}

\item[{\rm{(iii)}}]
\begin{equation} \lb{3.21}
\{x_i, x_j\}=0
\end{equation}

\item[{\rm{(iv)}}]
\begin{equation} \lb{3.22}
\{x_i, y_j\} = (W^{-1})_{ji}
\end{equation}
\end{SL}
Moreover, $\{y_i, y_j\} =0$ for all $i,j$ if and only if $U\equiv 0$.

Conversely, if a symplectic manifold has a Poisson bracket obeying \eqref{3.21}
and \eqref{3.22}, then $\omega$ has the form \eqref{3.19}.
\end{proposition}

\begin{remarks} 1. In \eqref{3.20}, $\omega^\ell$ means the $\ell$-fold wedge product
$\omega\wedge \dots \wedge \omega$.

\smallskip
2. $W$ and $U$ can be functions of $x,y$.

\smallskip 
3. These calculations do not use $d\omega =0$ (or, equivalently, the Jacobi identity). 
\end{remarks}

\begin{proof} Write
\[
v=\sum_{k=1}^\ell \biggl( r_k\, \f{\partial}{\partial x_k} + t_k\,
\f{\partial}{\partial y_k}\biggr)
\]
If $\vec{t}\neq 0$, pick $\ti s= W^{-1} t$ and $\ti v=\sum_{k=1}^\ell \ti s_k
\f{\partial}{\partial x_k}$. If $\vec{t}=0$, $\vec{r}\neq 0$, and with $\ti s
= W^{-1} r$, we can take $\ti v=\sum_{k=1}^\ell \ti s_k \f{\partial}{\partial y_k}$.
Thus, invertibility implies nondegeneracy. If $W$ is not invertible, then
$\det(W)=0$ and the calculation of (ii) shows $\omega$ is degenerate.

(ii) is an easy calculation: using the distributive law to expand $\omega^\ell$,
any term with a $dx_i\wedge dx_j$ term has at least $\ell+1$ $dx$'s, and so is zero
by antisymmetry. The only products of $\ell$ $dx_i\wedge dy_j$ that are nonzero
are of the form $[dx_{\pi(1)}\wedge dy_{\sigma(1)}] \wedge\dots\wedge
[dx_{\pi(\ell)}\wedge dy_{\sigma(\ell)}]$ with $\pi,\sigma$ permutations,
and this is  $(-1)^{\sigma\pi^{-1}} dx_1 \wedge dy_1 \wedge \dots \wedge
dy_\ell$. The sum over $\sigma$ for fixed $\pi$ yields $\det(W)$ and then
the sum over $\pi$ gives $\ell$!.

(iii), (iv). We begin by noting that if $\{\zeta_j\}_{j=1}^{2\ell}$ is any
coordinate system, if $\Omega$ is a $2\ell\times 2\ell$ antisymmetric real
matrix, and
\begin{equation} \lb{3.33}
\omega =\sum_{j,k=1}^{2\ell} \Omega_{jk} d\zeta_j\wedge d\zeta_k
\end{equation}
is a symplectic form, then for
\begin{equation} \lb{3.34}
H_{\zeta_j} =\sum_{k=1}^\ell \alpha_k \f{\partial}{\partial\zeta_k}
\end{equation}
we have
\[
\omega \biggl(H_{\zeta_j}, \f{\partial}{\partial\zeta_k}\biggr) = \delta_{jk}
\Rightarrow \sum_{m=1}^\ell \Omega_{mk} \alpha_m = \delta_{kj}
\]
or (since $\Omega_{mk}=-\Omega_{km}$)
\begin{equation} \lb{3.35}
\alpha_k = [(-\Omega)^{-1}\delta_j]_k = -(\Omega^{-1})_{kj}
\end{equation}
so that
\begin{equation} \lb{3.36}
\{\zeta_j, \zeta_k\} = H_{\zeta_j} (\zeta_k) = -(\Omega^{-1})_{kj}
\end{equation}

For our case if $\zeta = (x_1, \dots, x_\ell, y_1, \dots, y_\ell)$,
$\Omega$ has the $\ell\times \ell$ block form
\begin{equation} \lb{3.37}
\Omega = \left(
\renewcommand{\arraystretch}{1.5}
\begin{tabular}{c | c}
$U$ & $\,\,W\,\, \,$ \\
\hline
$-W^t$ & $0$
\end{tabular}
\right)
\end{equation}
for which
\begin{equation} \lb{3.38}
\Omega^{-1}= \left(
\renewcommand{\arraystretch}{1.5}
\begin{tabular}{c | c}
$ 0$ & $-W_t^{-1}$ \\
\hline
$\, W^{-1} \, $ & $W^{-1} UW_t^{-1}$
\end{tabular}
\right)
\end{equation}
so \eqref{3.21} and \eqref{3.22} follow from \eqref{3.36}.

For $\{y_i, y_j\}=0$ if and only if $U\equiv 0$, we note $\{y_i, y_j\}=0$ if
and only if $W^{-1} UW^{-1}=0$ by \eqref{3.18}.

For the converse, we note that \eqref{3.21} and \eqref{3.22} plus antisymmetry imply
$\Omega^{-1}$ has the form \eqref{3.38} for some $U$, and then \eqref{3.33} implies
\eqref{3.19}.
\end{proof}

From this proposition, we can find the basic symplectic two-form in both $(a,b)$
and $(\mu,x)$ coordinates:

\begin{theorem}\lb{T3.6} The symplectic form defined  by \eqref{1.3}/\eqref{1.4} has
the form:
\begin{equation} \lb{3.29}
4\sum_{1\leq \ell\leq k\leq N-1} a_\ell^{-1} \, da_\ell \wedge db_k
\end{equation}
\end{theorem}

\begin{remark} \eqref{1.3}/\eqref{1.5} only define a two-form. That it is symplectic 
(i.e., $d\omega =0$) follows from \eqref{3.29}. 
\end{remark}

\begin{proof} Since \eqref{1.3}/\eqref{1.4} say
\[
\{a_k, b_k\} = -\{a_k, b_{k+1}\} = \tfrac14\, a_k
\]
Proposition~\ref{P3.5} (which, as we noted, does not use that the Poisson bracket  
obeys the Jacobi identity) says that (recall that $b_1, \dots, b_{N-1}; a_1,
\dots, a_{N-1}$ are a set of coordinates)
\begin{equation} \lb{3.30x}
\omega =\sum_{i,j=1}^{N-1} W_{ij}\, da_i \wedge db_j
\end{equation}
where $(W^{-1})^t$ has the form
\begin{equation} \lb{3.30}
\begin{pmatrix}
\tfrac14\, a_1 & -\tfrac14\, a_1 & 0 & \dots \\
0 & \tfrac14\, a_2 & -\tfrac14\, a_2 & \dots \\
\dots & \dots & \dots & \dots
\end{pmatrix}
= D(1-N)
\end{equation}
where $D$ is the diagonal matrix $\f14\, a_k \delta_{k\ell}$ and $N= \left(
\begin{smallmatrix} 0 & 1 & 0 \\ 0 & 0 & 1 \\ \dots & \dots & \dots\end{smallmatrix}
\right)$ is the standard rank $N-2$ nilpotent on $\bbR^{N-1}$. Thus
\begin{align}
W &= [(D(1-N))^{-1}]^t = D^{-1} (1+N^t + (N^t)^2 + \cdots + (N^t)^{n-1}) \lb{3.31x} \\
&= \begin{pmatrix}
4a_1^{-1} & 0 & 0 & \dots \\
4a_1^{-1} & 4a_2^{-1} & 0 & \dots \\
4a_1^{-1} & 4a_2^{-1} & 4a_3^{-1} & \dots \\
\dots & \dots & \dots & \dots
\end{pmatrix} \lb{3.32x}
\end{align}
so \eqref{3.30x} is \eqref{3.29}.
\end{proof}

\begin{theorem}\lb{T3.7} The symplectic form defined by \eqref{1.3}/\eqref{1.4} has the
form:
\begin{equation} \lb{3.33x}
2\sum_{j=1}^N dx_j\wedge \rho_j^{-1}\, d\rho_j + \sum_{i,j} U_{ij} \, dx_i \wedge dx_j
\end{equation}
for some $U$.
\end{theorem}

\begin{remarks} 1. By $\sum_{j=1}^N x_j =c$ and $\sum_{j=1}^N \rho_j =1$, the $x$'s and
$\rho$'s are not independent, but they still define functions.

\smallskip
2. One could compute $U$ from \eqref{1.11a}.
\end{remarks}

\begin{proof} Define for $k=1, \dots, N-1$,
\begin{equation} \lb{3.34x}
y_k = \log \biggl[ \f{\rho_k}{\rho_N}\biggr]
\end{equation}
Then for $j=1, \dots, N-1$,
\begin{equation} \lb{3.35x}
\{x_j, y_k\} = \tfrac12\, \delta_{jk}
\end{equation}
since, by \eqref{1.11},
\begin{align*}
\{x_j,y_k\} &= \{x_j, \log \rho_k\} - \{x_j, \log \rho_n\} \\
&= (\tfrac12\, \delta_{jk} -\rho_j) - (-\rho_j)
\end{align*}

$\{x_j\}_{j=1}^{N-1}$, $\{y_j\}_{j=1}^{N-1}$ are local (indeed, global) coordinates
since $\sum_{j=1}^n \rho_j=1$ lets us invert \eqref{3.34x},
\begin{align}
\rho_j = \f{e^{y_j}}{[1+\sum_{\ell=1}^{N-1} e^{y_\ell}]}  \lb{3.36x} \\
\rho_N = \f{1}{[1+\sum_{\ell=1}^{N-1} e^{y_\ell}]} \lb{3.37x}
\end{align}
so $\{x_j, x_k\}=0$ and \eqref{3.36} plus Proposition~\ref{P3.5} imply
\begin{equation} \lb{3.38x}
\omega =2\sum_{j=1}^{N-1} dx_j \wedge dy_j + \sum_{i,j} U_{ij} dx_i \wedge dx_j
\end{equation}
\eqref{3.33x} follows if we note that
\[
dy_j = \rho_j^{-1} \, d\rho_j - \rho_N^{-1} \, d\rho_N
\]
and
\[
\sum_{j=1}^{N-1} dx_j \wedge d\rho_N = -dx_N \wedge d\rho_N
\]
since $\sum_{j=1}^N dx_j =0$.
\end{proof}

\section{Coefficient Stripping for Paraorthogonal Polynomials} \lb{s4}

Paraorthogonal polynomials are defined \cite{JNT89} by
\begin{equation} \lb{4.1}
P_N (z;\{\alpha_j\}_{j=0}^{N-2},\beta) = z\Phi_{N-1}(z) - \bar\beta
\Phi_{N-1}^*(z)
\end{equation}
where $\beta\in\partial\bbD$. Second kind paraorthogonal polynomials by
\begin{equation} \lb{4.2}
Q_N (z;\{\alpha_j\}_{j=0}^{N-2},\beta) =z\Psi_{N-1}(z) + \bar\beta
\Psi_{N-1}^*(z)
\end{equation}
where $\Psi$ are the second kind polynomials. These polynomials have been
extensively studied recently \cite{CMV,CMV06,Gol02,Sim308,Wong} and,
in particular, the second kind polynomials are introduced in \cite{Sim308,Wong}.
The symbols $P$ and $Q$ are new but quite natural. These are relevant for the
following reason:

\begin{proposition}\lb{P4.1} Let
\begin{equation} \lb{4.3}
d\mu = \sum_{j=1}^N \mu_j \delta_{z_j}
\end{equation}
be a pure point probability measure on $\partial\bbD$ with each $\mu_j >0$.
Let $\Phi_0,\Phi_1, \dots, \Phi_N$ be the monic OPUC where, in particular,
\begin{equation} \lb{4.4}
\Phi_N(z) = \prod_{j=1}^N (z-z_j)
\end{equation}
Then the recursion relations define parameters $\{\alpha_j\}_{j=0}^{N-2}$ and
$\alpha_{N-1} =\beta\in\partial\bbD$ where
\begin{equation} \lb{4.5}
\beta = (-1)^{N+1} \prod_{j=1}^N \bar z_j
\end{equation}
Here
\begin{equation} \lb{4.6}
\Phi_N (z,d\mu) = P_N (z;\{\alpha_j\}_{j=0}^{N-2},\beta)
\end{equation}

Moreover,
\begin{align}
F(z,d\mu) &\equiv \int \f{e^{i\theta}+z}{e^{i\theta}-z}\, d\mu(z) \notag \\
&=\sum_{j=1}^N \mu_j \, \f{z_j+z}{z_j-z} \lb{4.7} \\
&= - \f{Q_N(z; \{\alpha_j\}_{j=0}^{N-2},\beta)}{P_N(z; \{\alpha_j\}_{j=0}^{N-2},\beta)} \lb{4.8}
\end{align}
In particular,
\begin{equation} \lb{4.9}
Q_N(z) = \sum_{j=1}^N \mu_j (z+z_j) \prod_{\ell\neq j} (z-z_\ell)
\end{equation}
\end{proposition}

\begin{remarks} 1. Parts of this proof are close to results of Jones, Nj\r{a}stad, and Thron
\cite{JNT89} as presented in Theorem~2.2.12 of \cite{OPUC1}.

\smallskip
2. \eqref{4.9} is, of course, the analog of \eqref{3.5}.

\smallskip
3. It is known that $\Ima F(e^{i\theta})$ is strictly monotone and pure imaginary
between poles, so there is a single zero between poles. Thus, \eqref{4.8} shows
the zeros of $P_n$ and $Q_n$ interlace, a result of \cite{Sim308,Wong} obtained
by other means.
\end{remarks}

\begin{proof} If $d\mu$ has the form \eqref{4.3}, then $L^2 (\partial\bbD, d\mu)$ is
$N$-dimensional, so $\{\Phi_j\}_{j=0}^{N-1}$ span the whole space so $\Phi_N (e^{i\theta})=0$
for $d\mu$ a.e.\ $\theta$. It follows that $\Phi_N$ must obey (we use $P_N$ for later
purposes, even though we have not yet proven that $\Phi_N$ is a $P_N$)
\begin{equation} \lb{4.10}
P_N(z) = \prod_{j=1}^N (z-z_j)
\end{equation}

In particular, $\alpha_{N-1} = \ol{-\Phi_N(0)}$ is given by \eqref{4.5} and so is a
$\beta\in\partial\bbD$. Szeg\H{o} recursion thus says \eqref{4.1} holds.

Using $^*$ for degree $N$-polynomials, \eqref{4.1}/\eqref{4.2} imply
\begin{equation} \lb{4.11}
P_N^* = -\beta P_N \qquad Q_N^* =\beta Q_N
\end{equation}
so
\begin{equation} \lb{4.12}
-\f{Q_N}{P_N} = \f{Q_N^*}{P_N^*}
\end{equation}
and \eqref{4.8} is just (3.2.19) of \cite{OPUC1} (this is proven for $(\alpha_0, \dots,
\alpha_{N-1}, 0,0\dots)$ with $\abs{\alpha_{N-1}}<1$, but by taking $\alpha_{N-1}\to
\partial\bbD$, one gets $F=Q_N^*/P_N^*$). \eqref{4.7}, \eqref{4.8}, and \eqref{4.10}
imply \eqref{4.9}.
\end{proof}

Since we are about to consider changes in $N$, we shift notation from $N$ to $n$.
Since $P_n, Q_n$ are defined via a boundary condition on $\alpha_{n-1}$,
if we want to prove something inductively, we need to ``strip from the front,"
that is, remove $\alpha_0$. In \cite{OPUC1,OPUC2} we called this coefficient
stripping. For OPUC, this involves the formula for $\Phi_n (z;\{\alpha_j\}_{j=0}^{n-1})$
and $\Phi_n^* (z;\{\alpha_j\}_{j=0}^{n-1})$ in terms of $\alpha_0$ and $\Phi_{n-1}
(z;\{\alpha_{j+1}\}_{j=0}^{n-2})$ and $\Phi_{n-1}^*(z;\{\alpha_{j+1}\}_{j=0}^{n-2})$.
Quite remarkably (and conveniently for use in Sections~\ref{s5} and \ref{s6}),
coefficient stripping for $P_n, Q_n$ involves another $P,Q$ pair even though $P,Q$
are $(\Phi,\Psi)$ mixed! It turns out simpler to state things in terms of
\begin{align}
C_n (z;\{\alpha_j\}_{j=0}^{n-2},\beta) = \f{(P_n + Q_n)}{2} \lb{4.13} \\
S_n (z;\{\alpha_j\}_{j=0}^{n-2},\beta) = \f{(P_n - Q_n)}{2} \lb{4.14}
\end{align}
which we name analogously to $\cosh$ and $\sinh$.

\begin{theorem}\lb{T4.2} Let $\beta\in\partial\bbD$ be fixed and $\{\alpha_j\}_{j=0}^{n-2}
\in\bbD^{n-1}$. Let $C_n$ denote $C_n (z;\{\alpha_j\}_{j=0}^{n-2},\beta)$, and similarly for
$S_n$. Define $C_{n-1}$ by
\begin{equation} \lb{4.15}
C_{n-1} = C_{n-1} (z; \{\alpha_{j+1}\}_{j=0}^{n-3},\beta)
\end{equation}
{\rm{(}}with $\alpha_0$ removed{\rm{)}}, and similarly for $S_{n-1}$. Then
\begin{align}
C_n(z) &= z(C_{n-1}(z) - \alpha_0 z S_{n-1}(z)) \lb{4.16} \\
S_n(z) &=  - \bar\alpha_0 C_{n-1}(z) + S_{n-1}(z) \lb{4.17} \\
\end{align}
\end{theorem}

\begin{proof} As usual, let
\begin{equation} \lb{4.18}
\Lambda (\alpha,z) = \begin{pmatrix}
z & -\bar\alpha \\
-\alpha z & 1 \end{pmatrix}
\end{equation}
so with $\delta_+ = \binom{1}{0}\in\bbC^2$, $\delta_- =\binom{0}{1}\in\bbC^2$, Szeg\H{o}
recursion says
\begin{align}
\binom{\Phi_{n-1}}{\Phi_{n-1}^*} &= \Lambda (\alpha_{n-2}, z) \dots \Lambda
(\alpha_0, z) (\delta_+ + \delta_-) \lb{4.19} \\
\binom{\Psi_{n-1}}{-\Psi_{n-1}^*} &= \Lambda (\alpha_{n-2}, z) \dots
\Lambda (\alpha_0, z)(\delta_+ - \delta_-) \lb{4.20}
\end{align}

Let $B(\beta) = \left( \begin{smallmatrix} z & 0 \\ 0 & -\beta \end{smallmatrix}\right)$ and
conclude
\begin{equation} \lb{4.21x}
P_n = \langle (\delta_+ + \delta_-), (B(\beta) \Lambda_{n-2}(\alpha_{n-2},z) \dots \Lambda
(\alpha_0,z)) (\delta_+ + \delta_-)\rangle
\end{equation}
$Q_n$ is similar with the rightmost $(\delta_+ + \delta_-)$ replaced by $(\delta_+ - \delta_-)$.
Thus, $C_n,S_n$ correspond to using $\delta_+$ and $\delta_-$, and we have that
\begin{equation} \lb{4.22x}
\binom{C_n}{S_n} = \Lambda (\alpha_0,z)^t \Lambda (\alpha_1,z)^t \dots \Lambda
(\alpha_{n-1},z)^t B(\beta)^t \binom{1}{1}
\end{equation}
which implies
\begin{equation} \lb{4.21}
\binom{C_n}{S_n} = \Lambda (\alpha_0,z)^t \binom{C_{n-1}}{S_{n-1}}
\end{equation}
which is \eqref{4.16}--\eqref{4.17}.
\end{proof}

\begin{remark} One might think the transposes in \eqref{4.22x} should be adjoints,
but the shift of $\delta_\pm$ from the right to the left side of the inner product
introduces a complex conjugate.
\end{remark}

While we won't need it, we note the induced formula for $P_n,Q_n$:

\begin{theorem}\lb{T4.3} If $P_n(z) =P_n(z;\{\alpha_j\}_{j=0}^{n-2},\beta)$ and
$P_{n-1}(z) = P_n(z;\{\alpha_{j+1}\}_{j=0}^{n-3},\beta)$, and similarly for $Q_n,
Q_{n-1}$, then
\begin{align}
P_n &= \tfrac12 \, (z-\bar\alpha_0 + 1- \alpha_0 z) P_{n-1} +
\tfrac12\, (z-\bar\alpha_0 -1 + \alpha_0 z) Q_{n-1} \lb{4.22} \\
Q_n &= \tfrac12 \, (z+\bar\alpha_0 + 1+ \alpha_0 z) Q_{n-1} +
\tfrac12\, (z+\bar\alpha_0 -1 - \alpha_0 z) P_{n-1} \lb{4.23}
\end{align}
\end{theorem}

\begin{proof} This is immediate from \eqref{4.21}, together with
\begin{align*}
\binom{P_n}{Q_n} &=\left(
\begin{array}{rr} 1&1 \\ 1 & -1 \end{array}
\right) \binom{C_n}{S_n} \\
\binom{C_{n-1}}{S_{n-1}} &= \tfrac12\, \left( \begin{array}{rr} 1&1 \\ 1& -1
 \end{array}\right) \binom{P_{n-1}}{Q_{n-1}}
\end{align*}
and a calculation of
\[
\tfrac12\, \left(\begin{array}{rr} 1&1 \\ 1& -1 \end{array}\right) \Lambda^t
(\alpha_0,z)\left(\begin{array}{rr} 1&1 \\ 1& -1 \end{array}\right)
\qedhere
\]
\end{proof}

The algebraic simplicity of the recursion of $(C_n,S_n)$ relative to $(P_n,Q_n)$
is reminiscent of Schur vs.\ Carath\'eodory functions (see Section~3.4 of \cite{OPUC1}).
This is not a coincidence. Note that by \eqref{4.8}
\begin{align}
\f{C_n}{S_n} &= \f{P_n + Q_n}{P_n - Q_n} =
\f{1-[(-Q_n)/P_n]}{1+[(-Q_n)/P_n]} \lb{4.24} \\
&= \f{1-F}{1+F}\notag
\end{align}
so
\begin{equation} \lb{4.25}
f(z) = -\f{z^{-1} C_n}{S_n}
\end{equation}

Note that $S_n$ is a polynomial of degree $n-1$ and $C_n$ is a polynomial of degree $n$
vanishing at zero, so \eqref{4.25} is a ratio of polynomials of degree $n-1$. $-S_n$
and $z^{-1} C_n$ are thus essentially para-versions of the Wall polynomials. Indeed,
by the Pinter--Nevai formula \cite[Theorem~3.2.10]{OPUC1},
\begin{align*}
z^{-1} C_n(z) &= z B_{n-2}^*(z) + \bar\beta A_{n-2}(z) \\
S_n(z) &= -z A_{n-2}^*(z) - \bar\beta B_{n-2}(z)
\end{align*}

\section{Poisson Brackets of First and Second Kind OPUC} \lb{s5}

Our goal in this section is to prove \eqref{1.17} and \eqref{1.18}. Once we translate
to $C,S$, the induction will be even simpler than in the OPRL because $C_{n-1}, S_{n-1}$
are only dependent on $\{\alpha_{j+1}\}_{j=0}^{n-3}$ and so have zero Poisson bracket with
$\alpha_0$, so there are no analogs of terms like \eqref{2.1}. Since we are varying the index,
we use $n$ in place of $N$. We will begin by showing
\eqref{1.17}/\eqref{1.18} is equivalent to
\begin{align}
\{C_n (z), C_n(w)\} &= 0 = \{S_n(z), S_n(w)\} \lb{5.1} \\
\{C_n(z), S_n(w)\} &=-i \biggl( \f{C_n(z) wS_n(w) - C_n(w) z S_n(z)}{z-w}\biggr) \lb{5.2}
\end{align}

\begin{proposition} \lb{P5.1} \eqref{1.17} and the symmetry of $\{P_n(z), Q_n(w)\}$ under
$z\leftrightarrow w$ implies \eqref{5.1} and
\begin{equation} \lb{5.3}
\{C_n(z), S_n(w)\} = -\tfrac12\, \{P_n(z), Q_n(w)\}
\end{equation}
Similarly, \eqref{5.1} and the symmetry of $\{C_n(z), S_n(w)\}$ under $z\leftrightarrow w$
implies \eqref{1.17} and \eqref{5.3}.
\end{proposition}

\begin{proof} Given the antisymmetry of the Poisson bracket, the symmetry of $\{P_n(z),
Q_n(w)\}$ can be written
\[
\{P_n(z), Q_n(w)\} + \{Q_n(z), P_n(w)\} =0
\]
from which \eqref{1.17} implies \eqref{5.1}, and conversely, \eqref{5.1} and the symmetry
of $\{C_n(z), S_n(w)\}$ implies \eqref{1.17}.

The antisymmetry plus \eqref{1.17} also shows
\[
\{C_n(z), S_n(w)\} = \tfrac14\, 2\{Q_n(z), P_n(w)\}
\]
which is \eqref{5.3}, and the converse is similar.
\end{proof}

Next, we do a calculation showing the equality of the left side of \eqref{5.2} and
of \eqref{1.18} using $P_n = C_n + S_n$ and $Q_n = C_n -S_n$:
\begin{align}
& \text{RHS of \eqref{1.18}} \notag \\
&\, = -\f{i}{2}
\biggl[ \biggl( \f{z+w}{z-w}\biggr) [-2 C_n(z) S_n(w) + 2C_n(w) S_n(z)] \notag \\
& \qquad \qquad \qquad \qquad \qquad \qquad  + 2C_n(z) S_n(w) + 2C_n(w) S_n(z)\biggr] \notag \\
& \, = i \biggl[ \biggl\{\biggl( \f{z+w}{z-w}\biggr) -1 \biggr\} C_n(z) S_n(w)
 - \biggl\{\biggl( \f{z+w}{z-w}\biggr) + 1 \biggr\} C_n(w) S_n(z) \biggr] \notag \\
& \, = i \biggl[ \f{2C_n(z) wS_n(w) - 2C_n(w) zS_n(z)}{z-w}\biggr] \notag \\
& \, = -2 [\text{RHS of \eqref{5.2}}] \lb{5.3a}
\end{align}
consistent with \eqref{5.1}.

Thus, \eqref{1.17}/\eqref{1.18} is equivalent to \eqref{5.1}/\eqref{5.2}, which we proceed
to prove inductively, starting with $n=1$. We let $G_n(z,w)$ denote the right side of \eqref{5.2}.

\begin{proposition}\lb{P5.2} We have that
\begin{equation} \lb{5.4}
C_1(z) =z \qquad S_1(z) = -\bar\beta
\end{equation}
Thus, $G_1 (z,w)=0$ and \eqref{5.1}, \eqref{5.2} hold for $n=1$.
\end{proposition}

\begin{proof} Since
\begin{equation} \lb{5.5}
P_1(z) = z-\bar\beta \qquad Q_1(z) = z+\bar\beta
\end{equation}
\eqref{5.4} is immediate, and thus,
\[
G_1(z,w) = i\bar\beta \biggl( \f{zw-wz}{z-w}\biggr) =0
\]
Since $C_1, S_1$ are independent of $\{\alpha_j\}$, all Poisson brackets are zero and
\eqref{5.1}, \eqref{5.2} hold for $n=1$.
\end{proof}

\begin{remark} It is an interesting exercise to prove
\[
(P_1(z) Q_1(w) - P_1(w) Q_1(z)) \, \f{z+w}{z-w} - Q_1(z) Q_1(w) + P_1(z) P_1(w)
\]
is $0$ from \eqref{5.5}, and explain why $-QQ$ appears with a minus sign.
\end{remark}

\begin{lemma}\lb{L5.3} $G_n$ and $G_{n-1}$ are related by
\begin{equation} \lb{5.6}
G_n (z,w) = z\rho_0^2 G_{n-1}(z,w) - i \rho_0^2 zS_{n-1}(z) C_{n-1}(w)
\end{equation}
\end{lemma}

\begin{proof} By \eqref{4.16}, \eqref{4.17}, and the fact that
\begin{equation} \lb{5.7}
\calB (f,g)(z,w) = \f{f(z) g(w) - f(w) g(z)}{z-w}
\end{equation}
is zero if $f=g$:
\begin{align}
G_n(z,w) &= -i \calB(C_n, XS_n) (z,w) \notag \\
&= -i \calB (XC_{n-1} -\alpha_0 XS_{n-1}, -\bar\alpha_0 XC_{n-1} + XS_{n-1}) \notag \\
&= -i [\calB (XC_{n-1}, XS_{n-1}) + \abs{\alpha_0}^2 \calB(XS_{n-1}, XC_{n-1})] \notag \\
&= -i \rho_0^2 \calB (XC_{n-1}, XS_{n-1}) \lb{5.8}
\end{align}
where we use $Xf$ for the function
\[
(Xf)(z) = zf(z)
\]
and \eqref{5.8} using antisymmetry of $\calB$ in $f$ and $g$ and $1-\abs{\alpha_0}^2 =
\rho_0^2$.

Next, note
\begin{align*}
\calB(Xf, g) (z,w) &= \f{zf(z) g(w) -wf(w) g(z)}{z-w} \\
&= z \biggl( \f{f(z)g(w) -f(w) g(z)}{z-w}\biggr) + \f{(z-w)}{z-w}\, f(w) g(z) \\
&= z\calB (f,g) (z,w) + g(z) f(w)
\end{align*}
to conclude
\begin{equation} \lb{5.9}
-i \calB(XC_{n-1}, XS_{n-1}) = zG_{n-1}(z,w) -i zS_{n-1}(z) C_{n-1}(w)
\end{equation}
\eqref{5.8} and \eqref{5.9} imply \eqref{5.6}.
\end{proof}

\begin{theorem}\lb{T5.4} \eqref{5.1}--\eqref{5.2} hold for all $n$.
\eqref{1.17}--\eqref{1.18} hold for all $n$.
\end{theorem}

\begin{proof} As noted after Proposition~\ref{P5.1} and \eqref{5.3a}, \eqref{1.17},
\eqref{1.18} is equivalent to \eqref{5.1} and \eqref{5.2}. So, by Proposition~\ref{P5.2}
and induction, we need only prove
\begin{equation} \lb{5.10}
\text{\eqref{5.1}--\eqref{5.2} for } n-1 \Rightarrow \text{\eqref{5.1}--\eqref{5.2} for } n
\end{equation}

By \eqref{4.16} and \eqref{4.17}, $\{\alpha_0, \alpha_0\} = \{\bar\alpha_0, \bar\alpha_0\}=0$
the symmetry of $G_{n-1}(z,w)$ in $z\leftrightarrow w$, and the independence of $S_{n-1},
C_{n-1}$ on $\alpha_0, \bar\alpha_0$, we have
\[
\text{\eqref{5.1}--\eqref{5.2} for } n-1 \Rightarrow \text{\eqref{5.1} for } n
\]
(Notice how much simpler this is than the analog for OPRL!)

On the other hand, by \eqref{4.16} and \eqref{4.17},
\begin{align}
\{C_n(z), S_n(w)\} &= z\{C_{n-1}(z) -\alpha_0 S_{n-1}(z), - \bar\alpha_0 C_{n-1}(w) +
S_{n-1}(w)\} \notag \\
&= z\rho_0^2 G_{n-1} (z,w) + \{\alpha_0, \bar\alpha_0\} zS_{n-1}(z) C_{n-1}(w) \notag \\
&= z\rho_0^2 G_{n-1}(z,w) - i\rho_0^2 zS_{n-1}(z) C_{n-1}(w) \lb{5.11} \\
&= G_n (z,w) \lb{5.12}
\end{align}
\eqref{5.11} uses \eqref{1.16} and \eqref{5.12} follows from \eqref{5.6}. We have thus proven
\eqref{5.2} for $n$.
\end{proof}

\section{Fundamental Poisson Brackets for OPUC} \lb{s6}

In this section, we will establish \eqref{1.22}. Let $\{z_j\}_{j=1}^n$ be the zeros of
$P_n(z)$, $\theta_j$ given by \eqref{1.20} and $\mu_j$ by \eqref{1.19}. Unlike the
OPRL case where $x_1 < \cdots < x_n$ allowed global variables, these are only defined
locally, both because of the $z$'s not having a unique initial $z_1$ and because
of $2\pi$ ambiguities in $\theta_j$.

Throughout this section, we fix $\beta\in\partial\bbD$ and consider the $2n-2$-dimensional
manifold $\bbD^{n-1}$ of Verblunsky coefficients $\{\alpha_j\}_{j=0}^{n-2}$ and
$\alpha_{n-1}=\beta$. The associated finite CMV matrix $\calC$ has
\begin{equation} \lb{6.1a}
\det(\calC) =(-1)^{n-1} \beta
\end{equation}
The associated spectral measure given by \eqref{1.20} thus has
\begin{equation} \lb{6.1b}
\prod_{j=1}^n z_j = (-1)^{n-1} \beta
\end{equation}
or
\begin{equation} \lb{6.1cx}
\sum_{j=1}^n \theta_j = \arg ((-1)^{n-1}\beta)
\end{equation}
and, of course,
\begin{equation} \lb{6.1c}
\sum_{j=1}^n \mu_j =1
\end{equation}
so as local coordinates $\{\theta_j\}_{j=1}^{n-1}\cup \{\mu_j\}_{j=1}^{n-1}$ coordinatizes
$\bbD^{n-1}$. $\beta$ never appears explicitly below---this is not surprising since a rotation
can move $\beta$ to $1$.

As usual, one can replace $\Real\alpha, \Ima \alpha$ by ``independent" coordinates
$\alpha,\bar\alpha$ with
\begin{gather}
d\alpha = d(\Real\alpha) + i\, d(\Ima \alpha) \qquad d\bar\alpha = d(\Real \alpha) -
i\, d(\Ima\alpha)\lb{6.1d} \\
\f{\partial}{\partial\alpha} = \f12 \biggl[ \f{\partial}{\partial\Real\alpha} +
\f{1}{i}\, \f{\partial}{\partial\Ima\alpha}\biggr] \qquad
\f{\partial}{\partial\bar\alpha} = \f12 \biggl[ \f{\partial}{\partial\Real\alpha}
-\f{1}{i}\, \f{\partial}{\partial\Ima\alpha}\biggr] \lb{6.1e}
\end{gather}
Thus, the Poisson bracket defined by \eqref{1.16} can be defined in general by
\begin{equation} \lb{6.1f}
\{f,g\} = \sum_{j=1}^n i\, \rho_j^2 \biggl( \f{\partial f}{\partial\bar\alpha_j}\,
\f{\partial g}{\partial\alpha_j} - \f{\partial f}{\partial\alpha_j}\,
\f{\partial g}{\partial\bar\alpha_j}\biggr)
\end{equation}

We will denote
\begin{align}
d^2 \alpha_j &= d(\Real\alpha_j) \wedge d(\Ima\alpha_j) \lb{6.1g} \\
&= \f{i}{2}\, d\alpha_j\wedge d\bar\alpha_j \lb{6.1h}
\end{align}

\begin{proposition}\lb{P6.1} For OPUC,
\begin{equation} \lb{6.1}
\left.\{P_n(z), P_n(w)\}\right|_{z=z_k, \, w=z_j} = -z_k z_j
\{\theta_k,\theta_j\} \prod_{\ell\neq k} (z_k-z_\ell)
\prod_{\ell\neq j} (z_j-z_\ell)
\end{equation}
In particular,
\begin{equation} \lb{6.2}
\{\theta_k,\theta_j\}=0
\end{equation}
\end{proposition}

\begin{proof} Since $P_n(z) = \prod_{\ell=1}^n (z-z_\ell)$, we have
\[
\text{LHS of \eqref{6.1}} = \{z_k, z_j\} \prod_{\ell\neq k} (z_k-z_\ell)
\prod_{\ell\neq j} (z_j-z_\ell)
\]
\eqref{6.1} follows if we note
\[
\{e^{i\theta_k}, e^{i\theta_j}\} = (i)^2 e^{i\theta_k} e^{i\theta_j}
\{\theta_k,\theta_j\}
\]
\eqref{1.17} then implies \eqref{6.2}.
\end{proof}

\begin{proposition} \lb{P6.2} For OPUC,
\begin{equation} \lb{6.3}
\left.\{P_n(z), Q_n(w)\}\right|_{z=z_j, \, w=z_k} =
-i \{\theta_j,\mu_k\} 2z_j z_k \prod_{\ell\neq k} (z_j-z_\ell)
\prod_{\ell\neq j} (z_k-z_\ell)
\end{equation}
while
\begin{equation} \lb{6.4}
\begin{split}
 \lim_{\substack{z\to z_j \\ z\neq z_j}}\, &\left. \text{RHS of \eqref{1.18}}
\right|_{w=z_k} \\
& = -\f{i}{2} \{2z_k\mu_k\delta_{jk} (z_j + z_k) - 2z_j
2z_k \mu_j \mu_k\} \prod_{\ell\neq j} (z_j-z_\ell) \prod_{\ell\neq k}
(z_k-z_\ell)
\end{split}
\end{equation}
In particular,
\begin{equation} \lb{6.5}
\{\theta_j, \mu_k\} = \mu_j\delta_{jk} -\mu_j \mu_k
\end{equation}
\end{proposition}

\begin{proof} Given \eqref{4.9}, setting $z=z_j$, $w=z_k$ picks out the
$\{z_j,\mu_k\}$ term, that is,
\[
\text{LHS of \eqref{6.3}} = \{-z_j, \mu_k\} \biggl[\, \prod_{\ell\neq j} (z_j-z_\ell)\biggr]
\biggl[ 2z_k \prod_{\ell\neq k} (z_k-z_\ell)\biggr]
\]
where $2z_k$ is $\left. (w+z_k)\right|_{w=z_k}$. Since
\[
\{-z_j,\mu_k\} = -i\, z_j \{\theta_j,\mu_k\}
\]
we obtain \eqref{6.3}.

Setting $w=z_k$ in the right side of \eqref{1.18} gives (using $P(z_k)=0$, $Q(z_k) = \mu_k
(2z_k) \prod_{\ell\neq k} (z_k-z_\ell)$)
\[
\left.\text{RHS of \eqref{1.18}}\right|_{w=z_k} = -\f{i}{2} \biggl\{\f{P(z)}{z-z_k}
\, \mu_k (z+z_k) - Q(z) \mu_k\biggr\} (2z_k) \prod_{\ell\neq k} (z_k-z_\ell)
\]
Since $P(z)/(z-z_k) =\prod_{\ell\neq k} (z-z_\ell)$, the first term is zero as
$z\to z_j$ if $j\neq k$, and otherwise is $\prod_{\ell\neq j} (z_j -z_\ell)$.
This yields \eqref{6.4}. \eqref{6.4} plus \eqref{6.3} yields \eqref{6.5}.
\end{proof}

From $\{Q(z), Q(w)\}=0$, one can derive the formula for $\{\mu_j,\mu_k\}$. We do not
provide details.

Finally, for this section, we compute the symplectic form, $\omega$, in $d\alpha$, and
$d\theta_j,d\mu_j$ coordinates.

The $d\alpha$ coordinate is especially simple since in those coordinates, the Poisson
bracket is a product structure. By \eqref{1.16} and $\Real \alpha_j = \frac12
(\alpha_j + \bar\alpha_j)$, $\Ima\alpha_j = \f{1}{2i} (\alpha_j-\bar\alpha_j)$, we get
\[
\{\Real\alpha_j,\Ima\alpha_k\} = \f{\rho_j^2}{2}
\]
so by the $U=0$ case of Proposition~\ref{P3.5} (and \eqref{6.1h}):

\begin{proposition}\lb{P6.3} The symplectic form for OPUC with $\alpha_{n-1}=\beta\in
\partial\bbD$ fixed is
\begin{align}
\omega &= \sum_{j=0}^{n-2}\, \f{2}{\rho_j^2}\, d(\Real\alpha_j) \wedge
d(\Ima \alpha_j) \lb{6.6} \\
&=\sum_{j=0}^{n-2}\, \f{i}{\rho_j^2}\, (d\alpha_j \wedge d\bar\alpha_j) \lb{6.7}
\end{align}
\end{proposition}

The calculation for $\theta,\mu$ coordinates is essentially the same as for OPRL with
$x,\rho$ coordinates since $\sum_{j=1}^n \theta_j$ is constant. Since the $\f12$ in
\eqref{1.11} is missing from \eqref{1.22}, the $2$ in \eqref{3.33x} is absent and we get

\begin{proposition}\lb{P6.4} The symplectic form for OPUC in terms of $\{\mu_j, \theta_j\}$
is
\begin{equation} \lb{6.8}
\omega = \sum_{j=1}^N d\theta_j \wedge\mu_j^{-1}\, d\mu_j + \sum_{j,k} U_{jk}\,
d\theta_j \wedge d\theta_k
\end{equation}
for some $U$\!.
\end{proposition}

\section{Application 1: Jacobians of Coordinate Changes} \lb{s7}

This is the first of three sections in which we present applications of the fundamental
Poisson brackets \eqref{1.11} for OPRL and \eqref{1.22} for OPUC. Those in this section
and the next are not new, but we include them here because the literature on these
issues is not always so clear nor is it always emphasized that all that is needed are
the fundamental Poisson brackets.

Given two local coordinates $x_1, \dots, x_m$ and $y_1, \dots, y_m$, their Jacobian
$\abs{\f{\partial x}{\partial y}}$ is defined by
\begin{equation} \lb{7.1}
\biggl| \f{\partial (x)}{\partial(y)}\biggr| \equiv
\biggl| \det \biggl( \f{\partial x_j}{\partial y_k}\biggr)\biggr|
\end{equation}
named because
\begin{equation} \lb{7.2}
d^n x = \biggl| \f{\partial(x)}{\partial (y)}\biggr|\, d^n y
\end{equation}
relates the local volume elements. Of course, differential forms are ideal for this,
$dx_j =\sum \f{\partial x_j}{\partial y_k}\,dy_k$ means
\begin{equation} \lb{7.3}
dx_1\wedge \cdots \wedge dx_n = \det \biggl( \f{\partial x_j}{\partial y_k}\biggr)
dy_1 \wedge \dots \wedge dy_k
\end{equation}

Our goal here is to compute the two Jacobians $\abs{\f{\partial (a,b)}{\partial (x,\rho)}}$
for OPRL and $\abs{\f{\partial (\Real\alpha, \Ima\alpha)}{\partial (\theta,\mu)}}$ for
OPUC. It appears that these were first found for OPRL by Dumitriu--Edelman \cite{DE02}, and
motivated by that, for OPUC by Killip--Nenciu \cite{KN}. Their argument was indirect and
Killip--Nenciu asked if there weren't a direct calculation. This was provided by Forrester--Rains
\cite{FR} using forms and the continued fraction expansions of resolvents. Their argument
is not unrelated to the one below, but is lacking the connection to Poisson brackets and
we feel is more involved. In \cite{KN2}, Killip--Nenciu remark that Deift (unpublished)
explained to them how to go from the symplectic form for Jacobi parameters rewritten
in terms of $(x,\rho)$ to $\abs{\f{\partial(a,b)}{\partial (x,\rho)}}$ and they then
do the analogous OPUC calculation. As they remark, this is buried at the end of a long
paper, and the fact that only the form of Poisson brackets is involved is obscure.
We present this idea here to make the OPRL argument explicit and to present the OPUC
argument without a need for Lie algebra actions.

\begin{theorem}\lb{T7.1} On the $2N-2$-dimensional manifold of Jacobi parameters with
$\sum_{j=1}^N b_j$ fixed, use $\{a_j,b_j\}_{j=1}^{N-1}$ and $\{x_j,\rho_j\}_{j=1}^{N-1}$
as coordinates. Then
\begin{equation} \lb{7.4}
\biggl| \f{\partial (a,b)}{\partial (x,\rho)}\biggr| =
\biggl[\, \f{2^{-(N-1)}\prod_{j=1}^{N-1} a_j}{\prod_{j=1}^N \rho_j}\biggr]
\end{equation}
If one considers the $2N-1$-dimensional manifold without $\sum_{j=1}^N b_j$ fixed
{\rm{(}}equivalently, without $\sum_{j=1}^N x_j$ fixed{\rm{)}}, then \eqref{7.4}
is still true, but where the Jacobian is from $(a_1, \dots, a_{N-1}, b_1, \dots, b_N)$
to $(x_1, \dots, x_N, \rho_1, \dots, \rho_{N-1})$ rather than from $(a_1, \dots, a_{N-1},
b_1, \dots, b_{N-1})$ to $(x_j, \dots, x_{N-1}, \rho_j, \dots, \rho_{N-1})$.
\end{theorem}

\begin{proof} Let $\omega$ be the symmetric form on the $2N-2$-dimensional manifold determined
by \eqref{1.3}--\eqref{1.4}. By \eqref{3.29}, the $(N-1)$-fold wedge product
\begin{equation} \lb{7.5}
\omega\wedge\dots\wedge\omega = 4^{N-1} (N-1)! \prod_{j=1}^{N-1} a_j^{-1}\,
da_1\wedge db_1\wedge\dots\wedge da_{N-1}
\wedge db_{N-1}
\end{equation}

By \eqref{3.33x}, noting all terms with at least one $dx_k \wedge dx_j$ must have a repeated
$x$ and so vanish,
\begin{align}
\omega &\wedge  \dots\wedge\omega \notag \\
&= 2^{N-1} (N-1)! \sum_{j=1}^N \, \bigwedge_{k\neq j}
\rho_k^{-1} (dx_k\wedge d\rho_k) \notag \\
&= 2^{N-1} (N-1)! \biggl[\, \sum_{j=1}^N \biggl(\, \prod_{k\neq j} \rho_k^{-1}\biggr)\biggr]
 dx_1 \wedge d\rho_1 \wedge\dots\wedge dx_{N-1} \wedge d\rho_{N-1} \lb{7.6} \\
&= 2^{N-1}(N-1)! \biggl(\, \prod_{k=1}^N \rho_k^{-1}\biggr) dx_1 \wedge d\rho_1
\wedge \dots \wedge dx_{N-1} \wedge d\rho_{N-1} \lb{7.7}
\end{align}
Here \eqref{7.6} is obtained from
\[
dx_N =-\sum_{j\neq N} dx_j \qquad d\rho_N =-\sum_{j\neq N} d\rho_j
\]
and \eqref{7.7} uses
\begin{equation} \lb{7.8}
\sum_{j=1}^N \, \prod_{k\neq j} \rho_k^{-1} = \prod_{k=1}^N \rho_k^{-1}
\sum_{j=1}^N \rho_j =\prod_{k=1}^N \rho_k^{-1}
\end{equation}

\eqref{7.4} is immediate from the equalities of the right side of \eqref{7.5} and \eqref{7.7}.

To get the results for the $2N-1$-dimensional manifold, we note that for each fixed
$\sum_{j=1}^N x_j$, we have
\begin{equation} \lb{7.9}
\text{RHS of \eqref{7.7}} = \text{RHS of \eqref{7.5}}
\end{equation}
Since $\sum_{j=1}^N x_j = \sum_{j=1}^N b_j$, we can wedge the right side of \eqref{7.5}
with $db_1 + \cdots + db_N$ and of \eqref{7.7} by $dx_1 + \cdots + dx_N$ and so obtain the
analog of \eqref{7.9} but with the extra variable.
\end{proof}

\begin{theorem}\lb{T7.2} On the $2N-2$-dimensional manifold where $\alpha_{N-1}=\beta\in
\partial\bbD$ is fixed and
\[
\alpha_j = u_j + i\, v_j
\]
$(u,v\in\bbR)$, use $\{u_j,v_j\}_{j=0}^{N-2}$ and $\{\theta_j,\mu_j\}_{j=1}^{N-1}$ as
coordinates. Then
\begin{equation} \lb{7.10}
\biggl| \f{\partial (u,v)}{\partial (\theta,\mu)}\biggr| =
\f{2^{-(N-1)} \prod_{j=1}^N \rho_j^2}{\prod_{j=1}^N \mu_j}
\end{equation}
If we consider the $2N-2$-dimensional manifold without fixing $\beta$ and use
coordinates $\{u_j,v_j\}_{j=0}^{N-2}\cup\{\psi\}$ where $\beta = e^{i\psi}$
and $\{\theta_j\}_{j=1}^N \cup \{\mu_j\}_{j=1}^{N-1}$, then
\begin{equation} \lb{7.11}
\biggl| \f{\partial (u,v,\psi)}{\partial(\theta,\mu)}\biggr| = \text{RHS of \eqref{7.10}}
\end{equation}
\end{theorem}

\begin{proof} By \eqref{6.6}, the $N-1$-fold product
\begin{equation} \lb{7.12}
\omega\wedge\dots\wedge\omega = (N-1)! \, 2^{N-1} \biggl(\, \prod_{j=1}^{N-2}
\rho_j^{-2}\biggr) \, du_0\wedge dv_0 \wedge\dots\wedge du_{N-2}\wedge dv_{N-2}
\end{equation}
where, by \eqref{6.8} and the same calculation that led to \eqref{7.7} (i.e., \eqref{7.8}
with $\rho$ replaced by $\mu$),
\begin{equation} \lb{7.13x}
\omega\wedge\dots\wedge\omega = (N-1)! \prod_{j=1}^N \mu_j^{-1}\, d\theta_1\wedge d\mu_1
\wedge\dots\wedge d\theta_{N-1} \wedge d\mu_{N-1}
\end{equation}

\eqref{7.10} is immediate and \eqref{7.11} then follows from
\[
d\psi = \sum_{j-1}^N d\theta_j
\]
on account of \eqref{6.1c}.
\end{proof}

Note that in comparing \eqref{7.4} and \eqref{7.10} with Forrester--Rains \cite{FR} and
Dumitriu--Edelman \cite{DE02}, one needs to bear in mind that their $q_j$ are related
to $\rho_j$ (and $\mu_j$) by
\begin{equation} \lb{7.13}
\rho_j = q_j^2
\end{equation}
and thus
\begin{equation} \lb{7.14}
\f{d\rho_j}{\rho_j} = \f{2dq_j}{q_j}
\end{equation}
and therefore
\begin{equation} \lb{7.15}
\f{d\rho_1 \dots d\rho_{N-1}}{\rho_1 \dots \rho_N} = 2^{N-1} \, \f{1}{q_N}\,
\f{dq_1 \dots dq_{N-1}}{q_1\dots q_N}
\end{equation}
so they have no $2^{-(N-1)}$ in \eqref{7.4} but have a $2^{(N-1)}$ in \eqref{7.10}.
There is an extra factor of $1/q_N$ (which is in \cite{FR} but missing in \cite{DE02}
due to the fact that their $d^{N-1}q$ of \cite{DE02} is the measure on a sphere,
not Euclidean measure).

\section{Application 2: Generalized Toda and Schur Flows} \lb{s8}

In this section, we want to show how \eqref{1.11} allows the explicit solution of
flows like \eqref{1.6} and \eqref{1.7} in the $\rho$ (resp.\ $\mu$) variables.
As mentioned earlier, these results are not new. Theorem~\ref{T8.1} appears in
\cite{DLT85} and many times earlier. Theorem~\ref{T8.3} appears already in
\cite{KN2} as their Corollary~6.5.

\begin{theorem}\lb{T8.1} Let $f$ be a real-valued $C^1$ function on $\bbR$ and let
\begin{equation} \lb{8.1}
H=\sum_{j=1}^N f(x_j)
\end{equation}
and let $a,b$ solve {\rm{(}}with $\,\dot{} = d/dt${\rm{)}}
\begin{equation} \lb{8.2}
\dot{a}_k = \{H,a_k\} \qquad \dot{b}_k = \{H,b_k\}
\end{equation}
for some initial conditions for the Jacobi parameters. Then $\{a_k(t), b_k(t)\}_{k=1}^N$
are the Jacobi parameters associated to
\begin{equation} \lb{8.3}
x_j(t) =x_j(0) \qquad \rho_j(t) = \f{e^{\f12 tf'(x_k)}\rho_j(0)}
{\sum_{k=1}^N e^{\f12 tf'(x_k)} \rho_k(0)}
\end{equation}
\end{theorem}

\begin{remarks} 1. Since all that matters is $f$ at the $x_j$, we can restrict to
polynomials where if $f(x) =x^m$, then $H=\tr (J^m)$.

\smallskip
2. $\{x_j, x_k\}=0$ implies all $\tr (J^m)$ are constants of the motion.
\end{remarks}

\begin{example} \lb{E8.2} The Toda flow \eqref{1.6}/\eqref{1.7} corresponds to $H=2
\, \tr (J^2)$, so $f(x)=2x^2$ and $f'(x)=4x$, and \eqref{8.3} becomes
\begin{equation} \lb{8.4}
\rho_j(t) = \f{e^{2tx_j}\rho_j(0)}{\sum_{k=1}^N e^{2tx_k}\rho_k(0)}
\end{equation}
which is well known to solve Toda \cite{Mo75,Deift96}.
\end{example}

\begin{proof} The flow is generated by the vector field $\{H, \dott\}$ so
\begin{equation} \lb{8.5}
\dot{x}_j = \{H,x_j\} \qquad \dot{\rho}_j = \{H,\rho_j\}
\end{equation}
Since $\{x_j,x_k\}=0$, \eqref{8.5} implies $\dot{x}_j =0$, that is, $x_j(t) =x_j(0)$.

By \eqref{1.11}, if $j\neq N$,
\begin{align}
\biggl\{x_k, \f{\rho_j}{\rho_N}\biggr\} &= \rho_N^{-1} \{x_k, \rho_j\} -
\f{\rho_j}{\rho_N^2} \, \{x_j,\rho_N\} \notag \\
&= \tfrac12\, \rho_N^{-1} (\delta_{kj}\rho_j -\rho_j\rho_k) - \tfrac12\,
\rho_N^{-1} (\delta_{kN} \rho_j - \rho_j\rho_k) \notag \\
&= \tfrac12\, (\delta_{kj}-\delta_{kN}) \, \f{\rho_j}{\rho_N} \lb{8.6}
\end{align}
Thus
\begin{equation} \lb{8.7}
\f{d}{dt} \biggl( \f{\rho_j}{\rho_N}\biggr) = \tfrac12\, [f'(x_j) - f'(x_N)]
\f{\rho_j}{\rho_N}
\end{equation}

Defining $y_k(t)$ by \eqref{3.34x}, we get
\[
y_k (t) = y_k(0) + \tfrac12\, t (f'(x_j) -f'(x_N))
\]
which by \eqref{3.36x} yields \eqref{8.3}.
\end{proof}

\begin{theorem}\lb{T8.3} Let $f$ be a real-valued $C^1$ function on $\partial\bbD$
\begin{equation} \lb{8.8}
H=\sum_{j=1}^N f(e^{i\theta_j})
\end{equation}
and let $\alpha$ solve
\begin{equation} \lb{8.9}
\dot{\alpha}_j =\{H,\alpha_j\}
\end{equation}
with some initial conditions $\abs{\alpha_j(0)}<1$ and boundary conditions $\alpha_{-1}
=-1$, $\alpha_{N-1}=\beta$. Then $\{\alpha_j(t)\}_{j=1}^{N-2} \cup \{\alpha_{N-1}=\beta\}$
are the Verblunsky parameters associated to the $N$-point measure with parameters
\begin{equation} \lb{8.10}
\theta_j(t)=\theta_j(0) \qquad
\mu_j(t) = \f{e^{tg(e^{i\theta_j})}\mu_j(0)}
{\sum_{k=1}^N e^{tg(e^{i\theta_k})}\mu_k(0)}
\end{equation}
where $g$ is the function given by
\begin{equation} \lb{8.11}
g(e^{i\theta}) = \f{d}{d\theta}\, f(e^{i\theta})
\end{equation}
\end{theorem}

\begin{example} \lb{E8.4} If $f(e^{i\theta})=2\sin\theta$, then $H=\f1{i}
\tr (\calC-\calC^{-1})$ and \eqref{8.9} is the Schur flow \cite{AL75,AG94,FG99,Golppt,LTA}
\begin{equation} \lb{8.12}
\dot{\alpha}_j =\rho_j^2 (\alpha_{j+1} -\alpha_{j-1})
\end{equation}
$g$ is $2\cos\theta$ and \eqref{8.10} becomes
\begin{equation} \lb{8.13}
\mu_j(t) = \f{e^{2t\cos(\theta_j)} \mu_j(0)}
{\sum_{k=1}^N e^{2t\cos(\theta_k)}\mu_k(0)}
\end{equation}
\end{example}

\begin{proof} As above, $\{\theta_j,\theta_k\}=0$ implies $\theta_j(t) =\theta_j(0)$.
The calculation of $\mu_j(t)/\mu_N(t)$ is identical to the one in Theorem~\ref{T8.1}
except there is no factor of $\f12$ and $df(x)/dx$ is replaced by $df(e^{i\theta})/
d\theta$.
\end{proof}

By taking suitable limits, we get that certain measures on $\bbR$ (with all moments
finite) as initial conditions for difference equations on $\{a_j,b_j\}_{j=1}^\infty$
or $\{\alpha_j\}_{j=1}^\infty$ can be solved by
\[
\f{e^{tf'(x)}\, d\rho_{t=0}(x)}{\int e^{tf'(x)} \, d\rho_{t=0}(x)}
\]
See, for example, the discussion in \cite{LTA} and references therein.

\section{Differential Equations for OPRL} \lb{s9new}

We have just seen that symplectic flows are naturally defined on Jacobi parameters.
On the level of measures, they are solved by \eqref{8.3}. Here, we want to consider the
differential equations on OPs induced by these flows. In terms of \eqref{8.3}, we will suppose
$f$ is a polynomial; explicitly,
\begin{equation} \lb{n9.1}
\tfrac12\, f'(x) =\sum_{j=1}^k c_j x^j
\end{equation}
We drop the $c_0$ term since it drops out of $d\rho_t$. Here is our main result:

\begin{theorem}\lb{newT9.1} Under the map,
\begin{equation} \lb{n9.2}
d\rho_t (x) = \f{e^{\f12 f'(x)t} d\rho_0}{\int e^{\f12 f'(x)t }d\rho_0}
\end{equation}
with $f'$ given by \eqref{n9.1}, the monic orthogonal polynomials $P_n(x;d\rho_t)
=P_n(x;t)$ obey
\begin{equation} \lb{n9.3}
\dot{P}_n(x;t) = -\sum_{\ell=1}^k\, \sum_{j=1}^k c_j (J(t)^j)_{n+1,n-\ell+1}
\biggl(\, \prod_{j=n-\ell+1}^n a_j\biggr) P_{n-\ell}(x;t)
\end{equation}
where $J(t)$ is the Jacobi matrix of $d\rho_t$.
\end{theorem}

\begin{remark} The vector indices associated to $J$ start at $1$, so 
(with $p_j=P_j/\|P_j\|$) 
\begin{align}
J_{jk} &=\langle p_{j-1}, xp_{k-1}\rangle \lb{n9.3a} \\
J_{jj} &=b_j = \langle p_{j-1}, xp_{j-1}\rangle \lb{n9.3b} \\
J_{jj+1} &= J_{j+1\,\, j} =a_j = \langle p_j, xp_{j-1}\rangle \lb{n9.3c}
\end{align}
\end{remark}

\begin{example}\lb{newE9.2} Take $\f12 f'(x)=x$ (the Toda case). Then $J_{n+1, n}$ is all that
is relevant, and it is $a_n$, so \eqref{n9.3} becomes
\begin{equation} \lb{n9.4}
\dot{P}_n = -a_n^2 P_{n-1}
\end{equation}
This Toda case is known; see, for example, Peherstorfer \cite{Peh01}
and Barrios--Hern\'andez \cite{BHppt}.
\qed
\end{example}

\begin{example} \lb{new9.3} Take $\f12 f'(x)=x^2$. We have
\[
(J^2)_{n+1, n-1} = a_n a_{n-1} \qquad
(J^2)_{n+1,n} = a_n (b_{n+1} + b_n)
\]
so
\begin{equation} \lb{n9.5}
\dot{P}_n = -a_n^2 (b_{n+1} + b_n) P_{n-1} - a_n^2 a_{n-1}^2 P_{n-2}
\end{equation}
\qed
\end{example}

\begin{remark} Using $\langle P_n, P_n\rangle =a_na_{n-1}\dots a_1$ and \eqref{n9.3b},
one can easily go from differential equations for $P_n$ back to those for $a_n$ and
$b_n$.
\end{remark}

\begin{proof}[Proof of Theorem~\ref{newT9.1}] Since $P_n$ is monic, $\dot{P}_n$ is a
polynomial of degree $n-1$. Thus
\begin{align}
\dot{P}_n &=\sum_{\ell=1}^n \langle \dot{P}_n, p_{n-\ell}\rangle p_{n-\ell} \notag \\
&= -\sum_{\ell=1}^n \langle P_n, \tfrac12\, f' p_{n-\ell}\rangle p_{n-\ell} \lb{n9.6} \\
&= -\sum_{\ell=1}^n \langle p_n, \tfrac12\, f' p_{n-\ell}\rangle \f{\|P_n\|}{\|P_{n-\ell}\|} \,
P_{n-\ell} \notag \\
&= \text{RHS of \eqref{n9.3}} \notag
\end{align}
where \eqref{n9.6} comes from
\begin{align*}
O &= \f{d}{dt} \int P_n p_{n-\ell}\, d\mu_t  \\
&= \langle \dot{P}_n, p_{n-\ell}\rangle + \langle P_n, \dot{p}_{n-\ell}\rangle +
\langle P_n, \tfrac12\, f' p_{n-\ell}\rangle  \\
&= \langle \dot{P}_n, p_{n-\ell}\rangle + \langle P_n, \tfrac12\, f' p_{n-\ell}\rangle
\end{align*}
since $\dot{p}_{n-\ell}$ is a polynomial of degree $n-\ell$, and so orthogonal to $P_n$.
\end{proof}

One can also ask about derivatives of $p_n = P_n/\|P_n\|$. Since
\begin{equation} \lb{n9.7}
\dot{p}_n = \f{\dot{P}_n}{\|P_n\|} - \biggl[ \f{P_n}{\|P_n\|}\biggr]
\f{\|P_n\|^\mathbf{\cdot}}{\|P_n\|}
\end{equation}
we clearly need only find $\f{d}{dt}\log \|P_n\|$. Following \cite[Sec.~9.10]{OPUC2},

\begin{proposition}\lb{nP9.4} Under the map,
\begin{equation*}
d\rho_t (x) = \f{e^{\f12 f'(x)t} d\rho_0}{\int e^{\f12 f'(x)t}d\rho_0}
\end{equation*}
with $f'$ given by \eqref{n9.1}, the monic orthogonal polynomials $P_n(x;d\rho_t)
=P_n(x;t)$ obey
\begin{equation} \lb{n9.8}
\f{d}{dt}\, \log \|P_n\| = \tfrac12 \sum_{j=1}^k c_j
[(J^j)_{n+1\,\, n+1} - (J^j)_{1\, 1}]
\end{equation}
\end{proposition}

\begin{proof} Let
\begin{equation} \lb{n9.9}
N_t = \int e^{\f12 f'(x)t}\, d\rho_t
\end{equation}
Since $\dot{P}_n$ is orthogonal to $P_n$,
\begin{align}
\f{d}{dt}\, [\|P_n\|^2 N_t] &= \biggl( \int P_n^2 (\tfrac12\, f')\, d\rho_t\biggr) N_t \lb{n9.10} \\
&= \biggl[ \, \sum_{j=1}^k c_j (J^j)_{n+1\, \, n+1}\biggr] \|P_n\|^2 N_t \lb{n9.11}
\end{align}
since
\begin{equation} \lb{n9.12}
\int P_n^2 (\tfrac12\, f')\, d\rho_t = \|P_n\|^2
\langle p_n, (\tfrac12\, f') p_n\rangle
\end{equation}
and we have \eqref{n9.3a}. Similarly,
\begin{align}
\f{d}{dt}\, N_t &= N_t \int (\tfrac12\, f)\, d\rho_t \notag \\
&= \text{\eqref{n9.11} for } n=0 \lb{n9.13}
\end{align}

Since
\[
\f{d}{dt}\, \log \|P_n\| = \tfrac12\, \biggl[ \f{d}{dt}\, \log [\|P_n\|^2 N_t] -
\f{d}{dt}\, \log N_t\biggr]
\]
\eqref{n9.11} and \eqref{n9.13} imply \eqref{n9.8}.
\end{proof}

By \eqref{n9.3}, \eqref{n9.7}, and \eqref{n9.8}, we immediately have

\begin{theorem}\lb{nT9.5} Under the map,
\begin{equation*}
d\rho_t (x) = \f{e^{\f12 f'(x)t} d\rho_0}{\int e^{\f12 f'(x)t}d\rho_0}
\end{equation*}
with $f'$ given by \eqref{n9.1}, the orthonormal polynomials $p_n=P_n/\|P_n\|$ obey
\begin{equation} \lb{n9.14}
\dot{p}_n = -p_n \biggl( \tfrac12 \sum_{j=1}^k c_j [(J^j)_{n+1\,\, n+1} - (J^j)_{1\,\, 1}]\biggr)
-\sum_{\ell=1}^k \biggl[\, \sum_{j=1}^k c_j (J^j)_{n+1, n-\ell+1}\biggr] p_{n-\ell}
\end{equation}
\end{theorem}

\section{Differential Equations for OPUC} \lb{s10new}

For OPUC, one looks at flows generated by Laurent polynomials, $f(z)$, real on $\partial\bbD$.
If
\begin{equation} \lb{n10.1}
g(e^{i\theta}) = \f{\partial}{\partial\theta}\, f(e^{i\theta}) = \sum_{j=-\ell}^\ell
b_j e^{ij\theta}
\end{equation}
where $b_{-j}=\bar b_j$, then
\begin{equation} \lb{n10.2}
d\mu_t (\theta) = \f{e^{tg(e^{i\theta})}d\mu_0(\theta)}{\int e^{tg(e^{i\theta})}d\mu_0(\theta)}
\end{equation}
The analog of the $P_n$'s are the unnormalized CMV and alternate CMV bases
\begin{align}
Y_n &= \rho_0 \dots \rho_{n-1} \chi_n \lb{n10.3} \\
X_n &= \rho_0 \dots \rho_{n-1} x_n \lb{n10.4}
\end{align}
where $\{\chi_n\}_{n=0}^\infty$ and $\{x_n\}_{n=0}^\infty$ are the CMV and alternate
CMV bases (see \cite[Sec.~4.2]{OPUC1}). Then the exact same argument as led to
Theorem~\ref{newT9.1} leads to

\begin{theorem}\lb{newT10.1} Let $d\mu_t$ be given by \eqref{n10.2} with $g$ given by
\eqref{n10.1}. Then if $Y_n(t) = Y_n(d\mu_t)$,
\begin{equation} \lb{n10.5x}
\dot{Y}_n = -\sum_{m=0}^{n-1} \, \sum_{k=-\ell}^\ell b_k (\calC^k)_{mn}
\biggl(\, \prod_{j=m}^{n-1} \rho_j\biggr) Y_m
\end{equation}
with the same equation for $\dot{X}_n$ except the $(\calC^k)_{mn}$ is replaced by
$(\calC^k)_{nm}$.
\end{theorem}

\begin{remark} These are analogs of Proposition~\ref{nP9.4} and Theorem~\ref{nT9.5}.
\end{remark}

\begin{example}[Schur Flow]\lb{nE10.2} We have
\begin{equation} \lb{n10.5}
g(z)=z+ z^{-1}
\end{equation}
Ismail \cite{IsmBk} (see also Golinskii \cite{Golppt}) obtained the following differential
equation for the monic OPUC:
\begin{equation} \lb{n10.6}
\dot{\Phi}_n = \Phi_{n+1} - (z+\bar\alpha_n\alpha_{n-1}) \Phi_n -
\rho_{n-1}^2 \Phi_{n-1}
\end{equation}
We want to show this is equivalent to \eqref{n10.5x}. From our point of view, \eqref{n10.6}
is unusual since our arguments show that $\dot{\Phi}_n$ is a linear continuation of
$\{\Phi_j\}_{j=0}^{n-1}$ while \eqref{n10.6} has $\Phi_{n+1}$ and $\Phi_n$. (Since there
is a $z\Phi_n$, in fact, one can hope---and it happens---that RHS of \eqref{n10.6}
is $\sum_{j=0}^{n-1} c_j \Phi_j$.)

Since $\Phi_{n+1}=z\Phi_n -\bar\alpha_n \Phi_n^*$,
\begin{equation} \lb{n10.6x}
\text{RHS of \eqref{n10.6}}=-\bar\alpha_n (\Phi_n^* + \alpha_{n-1} \Phi_n)
-\rho_{n-1}^2 \Phi_{n-1}
\end{equation}
Since (see (1.5.41) of \cite{OPUC1}) $\Phi_n^* + \alpha_{n-1}\Phi_n =
\rho_{n-1}^2 \Phi_{n-1}^*$,
\begin{equation} \lb{n10.7}
\text{RHS of \eqref{n10.6}} = -\bar\alpha \rho_{n-1}^2 \Phi_{n-1}^* -
\rho_{n-1}^2 \Phi_{n-1}
\end{equation}
verifying that RHS of \eqref{n10.6} is indeed a polynomial of degree at most $n-1$.

\eqref{n10.7} mixes $Y$'s and $X$'s, so we use (see (1.5.40) of \cite{OPUC1})
$\Phi_{n-1}=\rho_{n-2}^2 z\Phi_{n-2} -\bar\alpha_{n-2} \Phi_{n-1}^*$ to find
\[
\text{RHS of \eqref{n10.6}} =-\rho_{n-1}^2 (\bar\alpha_n -\bar\alpha_{n-2})
\Phi_{n-1}^* -\rho_{n-1}^2 \rho_{n-2}^2 (z\Phi_{n-2})
\]which one can see is exactly \eqref{n10.6} or its $X_n$ analog (depending on
whether $n$ in $\Phi_n$ is even or odd).
\qed
\end{example}

\section{Application 3: Poisson Commutation for Periodic Jacobi and CMV Matrices} \lb{s9}

In this section, we want to show how the fundamental relations \eqref{1.11}/\eqref{1.22}
on finite Jacobi and CMV matrices yield a proof of the basic Poisson bracket relations
for the periodic case. For OPRL, the original proofs of this relation by Flaschka \cite{Fla1}
via direct calculation of Poisson brackets of eigenvalues is simple, but for OPUC, the
two existing proofs by Nenciu--Simon \cite{NS} and Nenciu \cite{NenDiss,Nen05} are
computationally involved.

Given $\{a_j,b_j\}_{j=1}^p\in\bbR^p$ or $\{\alpha_j\}_{j=0}^{p-1}\in\bbC^p$ (in the OPUC case,
we assume $p$ is even for reasons discussed in Chapter~11 of \cite{OPUC2}), we define
periodic Jacobi and CMV matrices by first extending the parameters periodically to $\bbZ$,
that is
\begin{equation} \lb{9.1}
a_{j+p} =a_j \qquad b_{j+p}=b_j \qquad \alpha_{j+p}=\alpha_j
\end{equation}
and then letting $J,\calC$ act on $\ell^2(\bbZ)$ as two-sided matrices with these
periodic parameters (two-sided CMV matrices are discussed on pp.~589--590 of \cite{OPUC2}).

In each case, $J,\calC$ act on $\ell^\infty(\bbZ)$ and commute with $\calS^p\colon
\{u_n\}\to \{u_{n+p}\}$. For each $e^{i\theta}\in\partial\bbD$, $J,\calC$ leave invariant
the $p$-dimensional space
\[
\ell_\theta^\infty = \{u\mid \calS^p u =e^{i\theta}u\}
\]
and so define $J(\theta),\calC(\theta)$ as operators on $\ell_\theta^\infty$. In terms of
the natural basis $\{\delta_j^\theta\}_{j=1}^N$,
\begin{equation} \lb{9.2}
(\delta_j^\theta)_k = \begin{cases}
e^{i\theta\ell} & k=p\ell +j \\
0 & k\not\equiv j \text{ mod } p \end{cases}
\end{equation}
$J(\theta)$ has a matrix of the form
\begin{equation} \lb{9.3}
\begin{pmatrix}
b_1 &  a_1 & 0 & \dots & e^{i\theta} a_p \\
a_1 & b_2 & a_2 & \dots & 0 \\
\vdots & \vdots & \vdots & {} & \vdots \\
e^{-i\theta} a_p & 0 & 0 & \dots & b_p
\end{pmatrix}
\end{equation}
of a Jacobi matrix with an extra term in the corners. The $\calC(\theta)$ are
described on pp.~719--720 of \cite{OPUC2}.

In each case, there is a natural transfer or monodromy matrix, $T_p(z)$ on $\bbC^2$,
so if $u\in\bbC^\infty$ solves $(J-z)u=0$ (resp.\ $(\calC-z)v=0$), then
\begin{equation} \lb{9.7}
T_p(z) \binom{u_1}{u_0} = \binom{u_{p+1}}{u_p}
\end{equation}
Constancy of Wronskians implies
\begin{equation} \lb{9.8}
\det (T_p(z)) = \begin{cases}
1 & \text{(OPRL)} \\
z^p & \text{(OPUC)} \end{cases}
\end{equation}
One defines the discriminant, $\Delta(z)$, by
\begin{equation} \lb{9.9}
\Delta(z) =\begin{cases}
\tr (T_p(z)) & \text{(OPRL)} \\
z^{-p/2} \tr (T_p(z)) & \text{(OPUC)}
\end{cases}
\end{equation}
(recall $p$ is even in the OPUC case). $z$ is clearly an eigenvalue of $J(\theta)$
(resp.\ $\calC(\theta)$) if and only if $e^{i\theta}$ is an eigenvalue of $T_p(z)$
(resp.\ $z^{-p/2} T_p(z)$; see \cite[p.~722]{OPUC2}). So, by \eqref{9.8},
\begin{equation} \lb{9.10}
z\text{ is an eigenvalue of $J(\theta)$ or $\calC(\theta)$} \Leftrightarrow
\Delta(z) =2\cos(\theta)
\end{equation}

For later purposes, we want to note that in neither case is $\Delta(z)$ automatically
monic. Rather,
\begin{equation} \lb{9.10a}
\Delta(z) = \begin{cases}
\f{1}{a_1\dots a_p}\, z^p + \text{lower order} & \text{(OPRL)} \\
\f{1}{\rho_1 \dots \rho_p}\, (z^{p/2} + z^{-p/2} + \text{orders in between}) & \text{(OPUC)}
\end{cases}
\end{equation}
(see \cite{Rice} and Chapter~11 of \cite{OPUC2}).

$\bbR^p \times (0,\infty)^{p-1}$ (resp.\ $\bbD^p$) are given Poisson brackets by
\begin{alignat}{2}
\{b_k, a_k\} &= -\tfrac14\, a_k \qquad && k=1,2,\dots, p \lb{9.4} \\
\{b_k,a_{k-1}\} &= \tfrac14\, a_{k-1} \qquad && k=2, \dots, p \lb{9.5} \\
\{b_1, a_p\} &= \tfrac14\, a_p \lb{9.6}
\end{alignat}
and
\[
\{\alpha_j,\alpha_k\}=\{\bar\alpha_j, \bar\alpha_k\} =0 \qquad
\{\alpha_j,\bar\alpha_k\} = -i\, \rho_j^2\, \delta_{jk} \qquad
j,k=0,\dots, p-1
\]

One big difference between OPRL and OPUC is that in the OPUC case, the center of the
Poisson bracket (i.e., the functions that Poisson commute with all functions)
is trivial, that is, $\{\cdot,\cdot\}$ defines a symplectic form.
But in the OPRL case,
\begin{equation} \lb{9.11}
\sum_{j=1}^p b_j \quad\text{and}\quad \prod_{j=1}^p a_j
\end{equation}
are easily seen to have zero Poisson brackets with all $a$'s and $b$'s, so the two-form
defined by $\{\cdot,\cdot\}$ is degenerate, but on the subspace where
\begin{equation} \lb{9.12}
\sum_{j=1}^p b_j =\beta \quad\text{and}\quad \prod_{j=1}^p a_j =\alpha
\end{equation}
the form is nondegenerate and the Poisson bracket defines a symplectic form.

Our main results in this section are:

\begin{theorem}[Flaschka \cite{Fla1}]\lb{T9.1} For OPRL:
\begin{equation} \lb{9.13}
\{\Delta(x),\Delta(y)\}=0
\end{equation}
If $\lambda_j(\theta)$ are the simple eigenvalues of $J(\theta)$, then
\begin{equation} \lb{9.14}
\{\lambda_j(\theta), \lambda_k(\theta')\}=0
\end{equation}
\end{theorem}

\begin{remarks} 1. \eqref{9.14} is normally only stated for $\theta=\theta'$, but
we will see it is immediate from \eqref{9.13} and \eqref{9.10}.

\smallskip
2. For $\theta\neq 0,\pi$, all eigenvalues are simple and $\lambda_j(\theta)$ are
global nonsingular functions. For $\theta=0,\pi$, there are subvarieties of codimension
$3$ where some eigenvalues are degenerate. $\lambda_j(\theta)$ are smooth away
from this subvariety and \eqref{9.14} only holds there.
\end{remarks}

\begin{theorem}[Nenciu--Simon \cite{NS}] \lb{T9.2} For OPUC:
\begin{equation} \lb{9.15}
\{\Delta(z),\Delta(w)\}=0
\end{equation}
If $\lambda_j(\theta)$ are the simple eigenvalues of $\calC(\theta)$, then
\begin{equation} \lb{9.16}
\{\lambda_j(\theta), \lambda_k(\theta')\} =0
\end{equation}
and
\begin{equation} \lb{9.17}
\biggl\{\, \prod_{k=0}^{p-1} \rho_k, \lambda_j(\theta)\biggr\} =0
\end{equation}
\end{theorem}

\begin{remarks} 1. As in the OPRL case, $\lambda_j(\theta)$ are simple if $\theta\neq 0,\pi$
and off a closed subvariety of codimension $3$ if $\theta=0$ or $\theta=\pi$.

\smallskip
2. $\prod_{k=0}^{p-1} \rho_k$ is the inverse leading term in $\Delta(z)$ for OPUC.
The analog for OPRL is $\prod_{k=1}^p a_k$. Of course, $\{\prod_{k=1}^p a_k,
\lambda_j(\theta)\}=0$, but we don't say it explicitly since for OPRL, $\prod_{k=1}^p a_k$
Poisson commutes with any function!
\end{remarks}

We note two general related aspects of these theorems. First, there are no variables
analogous to the $\rho$'s and $\mu$'s allowing a simple exact solution. The natural
complementary variables are the Dirichlet data which are defined specially on each
isospectral set rather than as a local function. The proper angle variables for this
situation are discussed in terms of the theory of Abelian integrals on a suitable
hyperelliptic surface \cite{DT76,DMN,GGH,KvM1,KvM2,KvM3,Kri,Mum,vM,vMM}.

Second, there is a sense in which the periodic situation is more subtle and interesting
than the finite $N$ structure. For finite $N$, the isospectral manifolds are open
simplexes parametrized by $\{\rho_j\}_{j=1}^N$ (or $\mu_j$) with $\sum_{j=1}^N \rho_j$=1.
Topologically this is $\bbR^{N-1}$. For the periodic case, the isospectral manifolds
are compact because $\prod_{j=1}^p a_j$ or $\prod_{j=1}^p \rho_j$ const.\ prevents
$a_j$ or $\rho_j$ from going to zero (for OPRL one also needs boundedness of 
$\sum a_j^2$). By the general theory of completely integrable
systems, these isospectral sets are tori generically of dimension one-half the
dimension of the symplectic manifold ($p-1$ for OPRL, $p$ for OPUC).

The idea of the proofs of these theorems is to cut off the periodic matrix to a finite
matrix and take limits of Poisson brackets. The limits of finite traces (normalized)
will be the moments of the density of states, not the periodic eigenvalues---but we
will see they are related. It will be useful to go back and forth between the two sets
of fundamental symmetric functions:
\begin{align}
t_k^{(N)} (\lambda_1, \dots, \lambda_N) &=\sum_{j=1}^N \lambda_j^k  \lb{9.18} \\
s_k^{(N)} (\lambda_1, \dots, \lambda_N) &=\sum_{i_1 < \cdots < i_k} \lambda_{i_1}
\dots \lambda_{i_k} \qquad 1\leq k\leq N \lb{9.19}
\end{align}
where we will need the following well-known combinatorial result:

\begin{proposition}\lb{P9.3} For $k\leq N$, we have
\begin{equation} \lb{9.20}
s_k^{(N)} (\lambda_1, \dots, \lambda_N)=\f{1}{k}\, t_k^{(N)} (\lambda_1, \dots,
\lambda_N) +r
\end{equation}
where $r$ is a polynomial in $\{t_j^{(N)}\}_{j=1}^{k-1}$ and also a polynomial in
$\{s_j^{(N)}\}_{j=1}^{k-1}$.
\end{proposition}

\begin{proof} Here is a simple proof. Let $A$ be the diagonal matrix with eigenvalues
$\lambda_1, \dots, \lambda_N$. Let $\wedge^k (\bbC^N)$ be the $k$-fold
antisymmetric product and $\otimes^k (\bbC^N)$ the tensor product (see, e.g., Appendix~A
of \cite{Groups}). Let $\pi\in\Sigma_k$, the symmetric group in $k$-objects, act as
$\vee_\pi$ on $\otimes^k (\bbC^N)$ via $\vee_\pi (x_1\otimes\cdots\otimes x_k)=
x_{\pi (1)}\otimes\cdots\otimes x_{\pi (k)}$. Then
\[
P_k = \f{1}{k!}\, \sum_{\pi\in \Sigma_k} (-1)^\pi \vee_\pi
\]
is the projection of $\otimes^k \bbC^N$ to $\wedge^k (\bbC^N)$. Thus,
\begin{align}
s_k^{(N)}(\lambda_1, \dots, \lambda_N) &= \tr_{\wedge^k (\bbC^N)} (\wedge^k (A)) \notag \\
&= \tr_{\otimes^k (\bbC^N)} (P_k\otimes^k A) \notag \\
&= \f{1}{k!} \, \sum_{\pi\in \Sigma_k} (-1)^\pi \tr (\vee_\pi \otimes^k A) \lb{9.21}
\end{align}
It is easy to see \cite[Appendix~A]{Groups} that if $\pi$'s disjoint cycle
decomposition has $\ell_1$ one-cycles, $\ell_2$ two-cycles, etc., then
\begin{align}
\tr (\vee_\pi \otimes^k A) &= \tr(A)^{\ell_1} \tr (A^2)^{\ell_2} \dots
\tr (A^k)^{\ell_k} \notag \\
&= \prod_{j=1}^k t_j^{(N)} (\lambda_1, \dots, \lambda_N)^{\ell_j}  \lb{9.22}
\end{align}
Since there are $(k-1)!$ permutations with one $k$-cycle, \eqref{9.21} and \eqref{9.22}
imply \eqref{9.20} with $r$ a polynomial in $t$'s. An easy induction then shows any $t_k$
is a polynomial in $\{s_j\}_{j=1}^k$ and so $r$ is also a polynomial in the $s$'s.
\end{proof}

The density of states of a periodic Jacobi or CMV matrix can be defined as the unique
measure, $d\gamma$, (on $\bbR$ or $\partial\bbD$) so that for all $k\in \{0,1,2,\dots\}$,
\begin{equation} \lb{9.23}
\int \lambda^k \, d\gamma(\lambda) = \int_0^{2\pi} \biggl[\, \f{1}{p} \sum_{j=1}^p
\lambda_j (\theta)^k \biggr] \f{d\theta}{2\pi}
\end{equation}
Moreover, it is easy to see \cite{Rice} that if $J_m (\{a_k\}_{k=1}^p, \{b_k\}_{k=1}^p)$
is the cutoff Jacobi matrix obtained by taking the two-sided infinite Jacobi matrix and
projecting onto the span of $\{\delta_j\}_{j=-m}^m$, then for any $k\in \{0,1,2,\dots\}$,
\begin{equation} \lb{9.24}
\f{1}{2m+1}\, \tr (J_m^k) \to \int \lambda^k \, d\gamma (\lambda)
\end{equation}
Similarly in the CMV case \cite{OPUC2}, if $\calC_m$ is defined by taking $\alpha_{-m}
\to -1$ and $\alpha_m\to\beta$ (for any $\beta)$, then
\begin{equation} \lb{9.25}
\f{1}{2m+1}\, \tr(\calC_m^k)\to \int \lambda^k\, d\gamma(\lambda)
\end{equation}

The key to the proof of Theorem~\ref{T9.1} is the following. Define for $\{a_j,
b_j\}_{j=1}^p$ fixed,
\begin{align}
t_k(\theta) &= \tr (J(\theta)^k) \lb{9.26} \\
s_k(\theta) &=\sum_{j_1< \cdots < j_k} \lambda_{j_1} (\theta) \dots
\lambda_{j_k}(\theta) \qquad k=0,1,\dots, p \lb{9.27x}
\end{align}.

\begin{theorem}\lb{T9.4} Consider OPRL of period $p$.
\begin{SL}
\item[{\rm{(i)}}] For $k\leq p-1$, $s_k(\theta)$ is independent of $\theta$.
\item[{\rm{(ii)}}] For $k\leq p-1$, $t_k(\theta)$ is independent of $\theta$.
\item[{\rm{(iii)}}]
\begin{equation} \lb{9.27}
s_p(\theta) = s_p(0) + (-1)^p \biggl(\, \prod_{j=1}^p a_j\biggr)
(2-2\cos\theta)
\end{equation}
\item[{\rm{(iv)}}]
\begin{equation} \lb{9.28}
t_p(\theta) = t_p(0) + (-1)^p p \biggl(\, \prod_{j=1}^p a_j\biggr) (2-2\cos\theta)
\end{equation}
\end{SL}
In particular, for $k\leq p$,
\begin{equation} \lb{9.29}
t_k(0) = \biggl[ (-1)^{p+1} (2p) \prod_{j=1}^p a_j\biggr]
\delta_{kp} + \int \lambda^k\, d\gamma (\lambda)
\end{equation}
\end{theorem}

\begin{remark} Directly from the form of the matrix \eqref{9.3}, one can see (ii) easily, and
(iv) with a little more work.
\end{remark}

\begin{proof} By \eqref{9.10a} and \eqref{9.10}, we see that
\begin{equation} \lb{9.30}
\prod_{j=1}^p (x-\lambda_j(\theta)) = \biggl(\, \prod_{j=1}^p a_j\biggr)
(\Delta(x) -2\cos\theta)
\end{equation}
which, expanding the products as
\begin{equation} \lb{9.31}
\prod_{j=1}^p (x-\lambda_j(\theta)) = x^p + \sum_{j=1}^p (-1)^j s_j(\theta)
x^{p-j}
\end{equation}
immediately implies (i) and (iii). (ii) and (iv) then follow from \eqref{9.20}, and
\eqref{9.29} then follows from \eqref{9.23}.
\end{proof}

For OPUC we have:

\begin{theorem}\lb{T9.5} Consider OPUC of even period $p$.
\begin{SL}
\item[{\rm{(i)}}] For $k=1,2,\dots, \f{p}{2}-1, \f{p}{2} +1, \dots, p-1$,
$s_k(\theta)$ is independent of $\theta$.
\item[{\rm{(ii)}}] For $k=1,2,\dots, \f{p}{2}-1$, $t_k(\theta)$ is independent of $\theta$.
\item[{\rm{(iii)}}]
\begin{equation} \lb{9.32}
s_{p/2}(\theta) = s_{p/2}(0) + (-1)^{p/2} \biggl(\,\prod_{j=0}^{p-1} \rho_j\biggr)
(2-2\cos\theta)
\end{equation}
\item[{\rm{(iv)}}]
\begin{equation} \lb{9.33}
t_{p/2}(\theta) = t_{p/2}(0) + (-1)^{p/2} \, \f{p}{2} \biggl(\,\prod_{j=0}^{p-1} \rho_j\biggr)
(2-2\cos\theta)
\end{equation}
\end{SL}
In particular, for $k\leq p/2$,
\begin{equation} \lb{9.34}
t_k(0) = (-1)^{p/2 +1} p \biggl(\,\prod_{j=0}^p \rho_j\biggr) \delta_{k \, p/2}
+ \int \lambda^k \, d\gamma(\lambda)
\end{equation}
\end{theorem}

\begin{remark} It will suffice to have control of $t_j$ and $\bar t_j$ for $0\leq \abs{j}
\leq p/2$ since $s_{p-j}=\bar s_j$, as we will see.
\end{remark}

\begin{proof} By \eqref{9.10a} and \eqref{9.10} and the fact that $z^{p/2}\Delta(z)$ is a
polynomial nonvanishing at $z=0$, we have that
\begin{equation} \lb{9.35}
\prod_{j=1}^p (z-\lambda_j(\theta)) = \biggl(\, \prod_{j=0}^{p-1}\rho_j\biggr)
(z^{p/2}\Delta(z) -2\cos(\theta) z^{p/2})
\end{equation}
which, via \eqref{9.31}, implies (i) and (iii). (ii) and (iv) then follow from \eqref{9.20},
and \eqref{9.34} from \eqref{9.23}.
\end{proof}

\begin{proof}[Proof of Theorem~\ref{T9.1}] Let
\begin{equation} \lb{9.37}
T_k^{(m)}(\{a_j,b_j\}_{j=0}^p) = \tr \bigl([J_m (\{a_k, b_k\}_{k=0}^p)]^k\bigr)
\end{equation}
and
\begin{equation} \lb{9.38}
\ti T_k^{(m)} (\{a_j,b_j\}_{j=-m}^m) = \tr \bigl( [J(\{a_k,b_k\}_{j=-m}^m)]^k\bigr)
\end{equation}
So $T_k^{(m)}$ is just $\ti T_k^{(m)}$ restricted to periodic sequences.

Let
\begin{equation} \lb{9.39}
t_k (\{a_j,b_j\}_{j=0}^p)=\int\lambda^k \, d\gamma(\lambda)
\end{equation}
As noted in \eqref{9.24}, $\f{1}{2m+1} T_k^{(m)}\to t_k$. This can be seen by noting
diagonal matrix elements of $(J_m)^k$ are uniformly bounded in $m$ for $k$ fixed and
equal to those of the infinite matrix so long as their index $j_0$ obeys $\abs{j_0}\leq
m-2k-1$.

This same equality and polynomial nature show for $1\leq j\leq p$,
\begin{equation} \lb{9.40}
\left. \f{\partial}{\partial a_{j+\ell p}} \, \ti T_k^{(m)}
\right|_{\substack{a_j = a_j^{(0)} \\ b_j=b_j^{(0)}}} =
\f{\partial}{\partial a_j}\, t_k (\{a_j^{(0)}, b_j^{(0)}\}_{j=1}^p
\end{equation}
where $\{a_j^{(0)},b_j^{(0)}\}$ is periodic, and similarly for $b$ derivatives.
\eqref{9.40} holds so long as as $\abs{j+\ell p} \leq m-2k-1$. It follows that
\begin{equation} \lb{9.41}
\left. \f{p}{m}\, \{\ti T_k^{(m)}, \ti T_j^{(m)}\}
\right|_{\substack{a=a^{(0)} \\ b=b^{(0)}}} \to \{t_k, t_\ell\}
\end{equation}

Since $\{\ti T_l^{(m)},\ti T_j^{(m)}\}=0$, we conclude
\begin{equation} \lb{9.42}
\{t_k,t_\ell\}=0
\end{equation}
for all $k,\ell$. By \eqref{9.29}
\begin{equation} \lb{9.42a}
\{t_k(0), t_\ell(0)\}=0
\end{equation}
for all $k,\ell$. Thus by \eqref{9.20}, we see for $1\leq k,\ell\leq p$,
\begin{equation} \lb{9.43}
\{s_k,s_\ell\}=0
\end{equation}

Since $\Delta(x)-2 = (a_1\dots a_p)^{-1} \prod_{j=1}^p (x-\lambda_j(0))$ and \eqref{9.31}
holds and $(a_1 \dots a_p)^{-1}$ is in the Poisson center, we see that
\begin{equation} \lb{9.44}
\{\Delta(x),\Delta(y)\}=0
\end{equation}
Thus
\begin{equation} \lb{9.45}
\{\Delta(x)-2\cos(\theta), \Delta(y) -2\cos (\theta')\}=0
\end{equation}
Evaluating this at $x=\lambda_j(\theta)$, $y=\lambda_k(\theta')$, we get \eqref{9.14}
as usual.
\end{proof}

\begin{proof}[Proof of Theorem~\ref{T9.2}] By repeating the arguments from the last proof,
we see that \eqref{9.42a} holds and that
\begin{equation} \lb{9.46}
\{t_k (0), t_\ell(0)\} = \{\ol{t_k(0)}, t_\ell(0)\}=0
\end{equation}

From this and Theorem~\ref{T9.5}, we see for $1\leq k,\ell\leq p/2$,
\begin{equation} \lb{9.47}
\{s_k(0), s_\ell(0)\} = \{s_k(0), \ol{s_\ell(0)}\} =0
\end{equation}
By \eqref{9.10a}, the product of the roots is $1$, that is,
\[
\prod_{j=1}^p \lambda_j(\theta) =1
\]
Since $\lambda_j\bar\lambda_j=1$ also, we see
\begin{equation} \lb{9.48}
\bar s_{p-j}(\theta) = \prod_{j=1}^p \lambda_j(\theta)\, \ol{s_{p-j}(\theta)} = s_j(\theta)
\end{equation}
so for $1\leq k,\ell \leq p-1$,
\begin{equation} \lb{9.49}
\{s_k(0), s_\ell(0)\}=0
\end{equation}

This, \eqref{9.35}, and \eqref{9.10a} imply that
\begin{equation} \lb{9.50}
\biggl\{\biggl(\, \prod_{j=0}^{p-1} \rho_j\biggr) \Delta(z), \biggl(\, \prod_{j=0}^{p-1}
\rho_j\biggr) \Delta(w)\biggr\} =0
\end{equation}

By a simple calculation if $\alpha_j=\abs{\alpha_j} e^{i\theta_j}$, then
\begin{equation} \lb{9.51}
\biggl\{\, \prod_{hj=0}^{p-1} \rho_j^2, g\biggr\} = -\prod_{j=0}^{p-1} \rho_j^2 \,
\sum_{j=0}^{p-1} \, \f{\partial g}{\partial \theta_j}
\end{equation}
so $-\log (\prod_{j=0}^{p-1}\rho_j^2)$ generates the Hamiltonian flow
\begin{equation} \lb{9.52}
\alpha_j\to \alpha_j e^{it}
\end{equation}
This transformation is implementable by a periodic unitary (see, e.g., \cite{CMVrev}), and
so it leaves the $\lambda_j(\theta)$ fixed. Thus, \eqref{9.17} holds, which implies
\begin{equation} \lb{9.53}
\biggl\{\,\prod_{j=0}^{p-1} \rho_j, \biggl(\, \prod_{j=0}^{p-1} \rho_j\biggr)
\Delta(z)\biggr\} =0
\end{equation}
so \eqref{9.50} implies \eqref{9.15}, which in turn implies \eqref{9.16}.
\end{proof}

\begin{remark} Unlike OPRL where $\prod_{j=1}^p a_j$ is in the Poisson center,
$\prod_{j=0}^{p-1}\rho_j$ generates a nontrivial flow, but one that leaves
$\lambda_j(\theta)$ invariant.
\end{remark}

\section{More Poisson Brackets for OPRL} \lb{s10}

In this section, we will discuss some additional Poisson brackets for OPRL and
related functions as a way of illuminating and extending the major results that
we proved earlier. We want to begin by showing that one can go backwards from
\eqref{1.11} to \eqref{1.14}/\eqref{1.15} (or, more precisely, to the Poisson
brackets \eqref{1.15} and $\{P_n(x), P_n(y)\}=0$; that $\{Q_n(x), Q_n(y)\}=0$
is then a consequence of \eqref{2.10}).

Clearly,
\begin{equation} \lb{10.1}
\{P_n(x), P_n(y)\}=0 \Leftrightarrow \{x_j,x_k\}=0
\end{equation}
To see that
\begin{equation} \lb{10.2}
\{x_j,\rho_k\}=\tfrac12\, [\delta_{jk}\rho_j - \rho_j\rho_k]\Leftrightarrow \text{\eqref{1.15}}
\end{equation}
we compute using \eqref{3.5} and $\{x_j, x_k\}=0$,
\begin{align}
2\{P_n(x), Q_n(y)\} &= 2\sum_{j,k} -\{x_j, \rho_k\}\,
\f{P_n(x)}{x-x_j}\, \f{P_n(y)}{y-y_k} \notag \\
&=\sum_{j,k} \rho_j\rho_k\, \f{P_n(x)}{x-x_j}\, \f{P_n(y)}{y-x_k} - \sum_j
\rho_j \, \f{P_n(x)P_n(y)}{(x-x_j)(y-x_j)} \lb{10.3}
\end{align}
The first sum gives $Q_n(x) Q_n(y)$ by \eqref{2.10} and the second sum is $(x-y)^{-1}
[P_n(x) Q_n(y) - P_n(y) Q_n(x)]$ if we note that
\begin{equation} \lb{10.4}
\f{1}{(x-x_j)}\, \f{1}{(y-x_j)} = \f{1}{x-y}\, \biggl[ \f{1}{y-x_j} - \f{1}{x-x_j}\biggr]
\end{equation}
This establishes \eqref{10.2}.

By \eqref{2.11} (with $n$ replaced by $n-1$) used in the first factor $\{P_n(x),
P_n(y)\}$, we get a relation among $\{P_{n-1}(x), P_n(y)\}$, $\{b_1, P_n(y)\}$,
$\{a_1^2, P_n(y)\}$, and $\{Q_{n-1}(x), P_n(y)\}$. Since the first three are computed
earlier, we get a formula for $\{Q_{n-1}(x), P_n(y)\}$, namely,
\begin{equation}\lb{12.4a}
\begin{split}
\{Q_{n-1}(x), &P_n(y)\} = 2P_{n-1}(x) - 2Q_{n-1}(x) - \\
& - b_1 a_1^{-2}
\biggl( \f{P_n(x) P_{n-1}(y) - P_{n-1}(x) P_n(y)}{x-y} - P_{n-1}(x) P_{n-1}(y)\biggr)
\end{split}
\end{equation}

Next, we want to discuss Poisson brackets with
\begin{equation} \lb{10.5}
m_n(x) = -\f{Q_n(x)}{P_n(x)}
\end{equation}
both to link to earlier work of Faybusovich--Gekhtman \cite{FG2000} and because one can
then take $n\to\infty$.

\begin{theorem}\lb{T10.1} We have that
\begin{align}
\{P_n(x), m_n(y)\} &= \tfrac12 \biggl[Q_n(x) m_n(y) -
\f{P_n(x) m_n(y) + Q_n(x)}{x-y}\biggr] \lb{10.6} \\
\{Q_n(x), m_n(y)\} &= \tfrac12 \biggl[ - Q_n(x) m_n(y)^2 + m_n(y)
\biggl[ \f{P_n(x) m_n(y) + Q_n(x)}{x-y}\biggr]\biggr] \lb{10.7}
\end{align}
\end{theorem}

\begin{proof} We have, since $\{P_n(x), P_n(y)\}=0$, that
\begin{equation} \lb{10.8}
\{P_n(x), m_n(y)\} = -P_n(y)^{-1} \{P_n(x), Q_n(y)\}
\end{equation}
The second term on the right of \eqref{1.15} leads to
\[
-\tfrac12\, P_n(y)^{-1} Q_n(x) Q_n (y) = \tfrac12\, Q_n(x) m_n(y)
\]
The first term on the right leads to
\[
\begin{split}
\tfrac12\, P_n(y)^{-1}  [P_n(x) Q_n(y) &- Q_n(x) P_n(y)](x-y)^{-1} \\
&\quad = (-\tfrac12\, P_n(x) m_n(y) - \tfrac12\, Q_n(x))(x-y)^{-1}
\end{split}
\]
proving \eqref{10.6}.

On the other hand, since $\{Q_n(x), Q_n(y)\}=0$,
\begin{align}
\{Q_n(x), m_n(y)\} &= Q_n(y) \{Q_n(x), -P_n(y)^{-1}\} \notag \\
&= m_n(y) P_n(y)^{-1} \{P_n(y), Q_n(x)\} \notag \\
&= m_n(y) P_n(y)^{-1} \{P_n(x), Q_n(y)\} \lb{10.9} \\
&= m_n(y) \{P_n(x), m_n(y)\} \lb{10.10}
\end{align}
so \eqref{10.6} implies \eqref{10.7}. \eqref{10.9} follows from the
symmetry of $\{P_n(x), Q_n(y)\}$ under $x\leftrightarrow y$, and \eqref{10.10}
from \eqref{10.8}.
\end{proof}

\begin{theorem}[\cite{FG2000}]\lb{T10.2} We have that
\begin{equation} \lb{10.11}
\{m_n(z), m_n(w)\} = \tfrac12\, (m_n(z) - m_n(w)) \biggl[ - \f{m_n(z)-m_n(w)}{z-w}
+ m_n(z) m_n(w)\biggr]
\end{equation}
\end{theorem}

\begin{proof} \ Since $\{P_n(z), P_n(w)\}=0=\{Q_n(z),Q_n(w)\}$ and  $\{P_n(z), Q_n(w)\}$ is symmetric
under $z\leftrightarrow w$ and $\{f,g\}$ is antisymmetric under $f\leftrightarrow g$,
\begin{align}
\{m_n(z), m_n(w)\} &= -\f{Q_n(z)}{P_n(w) Q_n(z)^2}\,
\{P_n(z), Q_n(w)\} - (z\leftrightarrow w) \lb{10.12} \\
&= \f{m_n(z)}{2} \biggl[ m_n(z) m_n(w) - \f{m_n(z)-m_n(w)}{z-w}\biggr] -
(z\leftrightarrow w) \notag
\end{align}
which is \eqref{10.11}.
\end{proof}

As $n\to\infty$, $m_n$ has a limit so long as $\{a_j, b_j\}_{j=0}^\infty$ are bounded
(actually so long as the moment problem is determinate), so \eqref{10.11} holds if
$m_n(z)$ is replaced by $m(z)$ for semi-infinite Jacobi matrices.

The unnormalized transfer matrix has the form
\[
\begin{pmatrix}
P_n & Q_n \\ P_{n-1} & Q_{n-1}
\end{pmatrix}
\]
Of the sixteen Poisson brackets for these four functions, we have computed twelve.
It would be interesting to know $\{Q_n(z), P_{n-1}(w)\}$.
Similarly, it would be interesting to know the Poisson brackets for the periodic
transfer matrix and use them to prove \eqref{9.13}.

\section{More Poisson Brackets for OPUC} \lb{s11}

We begin by showing that one can go backwards from \eqref{1.22} to \eqref{1.17}/\eqref{1.18}
(or, more precisely, to \eqref{1.18} and $\{P_n(z), P_n(w)\}=0$; $\{Q_n(z), Q_n(w)\}$ then
follows from the symmetry discussed after \eqref{1.18}).

Clearly,
\begin{equation} \lb{11.1}
\{P_n(z), P_n(w)\} = 0 \Leftrightarrow \{\theta_j,\theta_k\} =0
\end{equation}
To see that
\begin{equation} \lb{11.2}
\{\theta_j,\mu_k\} = \mu_j \delta_{jk}- \mu_j\mu_k\Rightarrow\text{\eqref{1.18}}
\end{equation}
we first need

\begin{lemma} \lb{L11.1} For any distinct, $z,w\in\bbC$ and $e^{i\theta}\in\bbD$,
\begin{equation} \lb{11.10a}
\f{e^{i\theta} +z}{e^{i\theta}-z}\, \f{e^{i\theta} +w}{e^{i\theta}-w} =
\f{z+w}{z-w} \biggl( \f{e^{i\theta} + z}{e^{i\theta} -z} -
\f{e^{i\theta} + w}{e^{i\theta} -w}\biggr) + 1
\end{equation}
\end{lemma}

\begin{proof} Consider first $e^{i\theta} =1$. Then
\[
\f{1+z}{1-z}\, \f{1+w}{1-w} -1 = \f{2(z+w)}{(1-z)(1-w)}
\]
and
\[
\f{1+z}{1-z} - \f{1+w}{1-w} = \f{2(z-w)}{(1-z)(1-w)}
\]
which implies \eqref{11.10a} for $e^{i\theta}=1$.

In that formula, replace $z$ by $e^{-i\theta} z$ and $w$ by $e^{-i\theta} w$
and thereby obtain the formula for general $e^{i\theta}$.
\end{proof}

Now we compute using \eqref{4.9} with $z_j=e^{i\theta_j}$,
\begin{align}
\{P_n(z), & Q_n(w)\} \notag \\
 &= \sum_{j,k} -\{z_j,\mu_k\}\, \f{P_n(z)}{z-z_j}\, P_n(w) \,
 \f{w+z_k}{w-z_k} \lb{11.3} \\
 &= -i \sum_{j,k} z_j [\delta_{jk}\mu_j-\mu_j\mu_k] \, \f{P_n(z)}{z-z_j} \,
 P_n(w)\, \f{w+z_k}{w-z_k} \notag \\
 &= -\f{i}{2} \sum_{j,k} [\delta_{jk} \mu_j -\mu_j\mu_k]
 \biggl( \f{e^{i\theta_j}+z}{e^{i\theta_j}-z}\biggr)
 \biggl( \f{e^{i\theta_k}+w}{e^{i\theta_k}-w}\biggr) P_n(z) P_n(w) \lb{11.4} \\
 &= \text{RHS of \eqref{1.18}} \lb{11.5}
\end{align}
verifying \eqref{11.2}. In the above, we get \eqref{11.4} by using
\begin{equation} \lb{11.6}
\f{z_j}{z-z_j} = \tfrac12 \biggl[ \f{z_j+z}{z_j-z} -1\biggr]
\end{equation}
and noting that the $-1$ term has $j$ dependence only from $[\delta_{jk}\mu_j -
\mu_j\mu_k]$ which sums to zero for each fixed $k$. To get \eqref{11.5}, we note
that \eqref{4.9} says
\begin{equation} \lb{11.7}
Q_n(z) = -\sum_k \mu_k  \biggl( \f{e^{i\theta_j}+z}{e^{i\theta_j}-z}\biggr) P_n(z)
\end{equation}
and use Lemma~\ref{L11.1} on the $\delta_{jk}\mu_j$ term. The first term in \eqref{11.10a}
gives $\f{z+w}{z-w} [P_n(z) Q_n(w) - Q_n(z) P_n(w)]$ by \eqref{11.7}. The second term in
\eqref{11.10a} gives $P_n(z) P_n(w)$ if we note that $\sum_{j,k} \mu_j \delta_{jk} =
\sum_j \mu_j =1$. The $\mu_j \mu_k$ term in \eqref{11.4} gives the $-Q_n(z) Q_n(w)$ by
\eqref{11.7}.

Next, we consider Poisson brackets with $F_n(z)$ given by \eqref{4.8} and $f_n(z) =
z^{-1} [\f{F_n(z)-1}{F_n(z)+1}]$.

\begin{theorem}\lb{T11.1} We have
\begin{align}
\{P_n(z), F_n(w)\} &= -\f{i}{2} \biggl[ (P_n(z) m_n(w) + Q_n(z)) \biggl(\f{z+w}{z-w}\biggr) \notag \\
& \qquad \qquad \qquad \qquad -P_n(z) -Q_n(z)m_n(w)\biggr] \lb{11.8} \\
\{Q_n(z), F_n(w)\} &= -F_n(w) (\text{RHS of \eqref{11.8}}) \lb{11.9} \\
\{P_n(z), f_n(w)\} &= (1-wf(w)) W_n(z,w) \lb{11.10} \\
\{Q_n(z), f_n(w)) &= -(1+wf(w)) W_n(z,w) \lb{11.11}
\end{align}
where
\begin{equation} \lb{11.12}
W_n(z,w) =-\f{i}{2} \biggl[ \f{(P_n(z) + Q_n(z)) - z(P_n(z) -Q_n(z)) f_n(w)}{z-w}\biggr]
\end{equation}
\end{theorem}

\begin{proof} \eqref{10.8} is valid with $m_n$ replaced by $F_n$. Thus
\begin{align}
\{P_n(z), F_n(w)\} & = -\f{i}{2} \biggl[ - P_n(w)^{-1} \biggl[ \biggl( \f{z+w}{z-w}\biggr)
(P_n(z) Q_n(w) - P_n(w) Q_n(z))\notag \\
& \qquad  + P_n(z) P_n(w) - Q_n(z) Q_n(w)\biggr]\biggr] \lb{11.13} \\
&= \text{RHS of \eqref{11.8}} \notag
\end{align}

\eqref{10.10} is also valid with $m_n$ replaced by $F_n$, so \eqref{11.9} follows from
\eqref{11.8}.

For the formula involving $f_n$, we use \eqref{4.25}. Analogous to \eqref{10.10}, one finds
\begin{equation} \lb{11.14}
\{C_n(z), f_n(w)\} = -wf(w) \{S_n(z), f(w)\}
\end{equation}
and, by \eqref{5.2},
\begin{align}
\{S_n(z), f_n(w)\} &= -w^{-1} S_n(w)^{-1} \{S_n(z), C_n(w)\} \lb{11.15} \\
&= -i \biggl[ \f{C_n(z) + zS_n(z) f_n(w)}{z-w}\biggr] \lb{11.16} \\
&= W(w,z) \lb{11.17}
\end{align}
by \eqref{4.13}/\eqref{4.14}.

We then get \eqref{11.10}, \eqref{11.11} using $P_n = C_n + S_n$, $Q_n = C_n-S_n$ and
\eqref{11.14}, \eqref{11.17}.
\end{proof}

\begin{theorem}[Gekhtman--Nenciu \cite{GNprep}] \lb{T11.2} We have that
\begin{align}
\{F_n(z), F_n(w)\} &= -\f{i}{2} \, (F_n(z) - F_n(w)) \biggl[ \biggl( \f{z+w}{z-w}\biggr) \notag \\
& \qquad \qquad (F_n(z) -F_n(w)) + 1 - F_n(z) F_n(w)\biggr] \lb{11.18} \\
\{f_n(z), f_n(w)\} & =-i\, \f{f(z)-f(w)}{z-w}\, (zf(z) -wf(w)) \lb{11.19}
\end{align}
\end{theorem}

\begin{remarks} 1. \cite{GNprep} has $-i$ and $-2i$ where we have $-\f{i}{2}$ and $-i$
because their $\{\cdot,\cdot\}$ is twice ours since, following \cite{KN2}, they dropped
the normalization of \cite{NS} which we keep.

\smallskip
2. While we will separately derive \eqref{11.18} and \eqref{11.19}, it is an illuminating
calculation to go from one to the other using $F(z) = (1+zf(z))/(1-zf(z))$.
\end{remarks}

\begin{proof} By the same calculation that led to \eqref{10.12},
\begin{align}
\{F_n(z), F_n(w)\} &= -\f{Q_n(z)}{P_n(w) P_n(z)^2}\, \{P_n(z),Q_n(w)\}
- (z\leftrightarrow w) \lb{11.20} \\
&= -\f{i}{2}\, F_n(z) \biggl[ \biggl( \f{z+w}{z-w}\biggr) (F_n(z) - F_n(w)) \notag \\
&\qquad \qquad\qquad +1 - F_n(z) F_n(w)\biggr] - (z\leftrightarrow w) \lb{11.21} \\
&= \text{RHS of \eqref{11.18}} \notag
\end{align}
where \eqref{11.21} follows from \eqref{1.18}.

Similarly, by \eqref{4.25} and the same calculation that led to \eqref{10.12},
\begin{align}
\{f_n(z), f_n(w)\} &= -(zw)^{-1} \, \f{C_n(z)}{S_n(z)^2 S_n(w)} \,
\{S_n(z), C_n(w)\} - (z\leftrightarrow w) \notag \\
&= -f_n(z) \, \f{\{C_n(z), S_n(w)\}}{wS_n(z) S_n(w)}  - (z\leftrightarrow w) \notag \\
&= \f{i\, f_n(z)}{z-w}\biggl[ \f{C_n(z)}{S_n(z)} - \f{z}{w}\, \f{C_n(w)}{S_n(w)}\biggr]
- (z\leftrightarrow w) \lb{11.22} \\
&= - \f{i\, f_n(z)}{z-w}\biggl[ zf_n(z) - \f{z}{w}\, wf_n(w)\biggr]
- (z\leftrightarrow w) \notag \\
&= -\f{i\, zf_n(z)}{z-w}\, (f_n(z) - f_n(w)) - (z\leftrightarrow w) \notag \\
&= \text{RHS of \eqref{11.19}} \notag
\qedhere
\end{align}
\end{proof}

Finally, we want to note that one can compute $\{\Phi_n(z), \Psi_n(z)\}$,
$\{\Phi_n^*(z), \Phi_n^*(w)\}$, and $\{\Phi_n(z), \Phi_n^*(w)$, $\{\Psi_n(z),
\Psi_n^*(w)\}$. We did not know $\{\Phi_n(z),\Psi_n^*(z)\}$ or $\{\Psi_n(z),
\Phi_n^*(z)\}$.  Nenciu \cite{Nenppt}, motivated by a preliminary version of 
this paper, has found this last missing bracket. With her formulae, one can get 
formulae for the brackets of the elements of the transfer matrix.

\begin{theorem}\lb{T11.3} We have
\begin{align}
\{\Phi_n(z), \Phi_n(w)\} &= \{\Psi_n(z), \Psi_n(w)\} = \{\Phi_n^*(z), \Phi_n^*(w)\} \notag \\
&= \{\Psi^*_n(z), \Psi^*_n(w)\} =0 \lb{11.23} \\
\{\Phi_n^*(z), \Psi_n^*(w)\} &= -\f{i}{2} \biggl[ (\Phi_n^*(z) \Psi_n^*(w) -
\Psi_n^*(z) \Phi_n^*(w)) \notag \\
&\qquad\quad  \biggl( \f{z+w}{z-w}\biggr) - \Phi_n^*(z) \Phi_n^*(w) +
\Psi_n^*(z) \Psi_n^*(w)\biggr] \lb{11.24} \\
\{\Phi_n(z), \Psi_n(w)\} &= -\f{i}{2}\biggl[ (\Phi_n(z) \Psi_n(w) - \Psi_n(z) \Phi_n(w)) \notag \\
& \qquad \quad \biggl( \f{z+w}{z-w}\biggr) + \Phi_n(z)\Phi_n(w) - \Psi_n(z) \Psi_n(w)\biggr] \lb{11.25}
\end{align}
\end{theorem}

\begin{proof} \eqref{1.17} and \eqref{1.18} depend on $\beta$ and are quadratic polynomials
in $\bar\beta$. Equality for all $\bar\beta\in\partial\bbD$ implies equality for all
$\bar\beta$, and so equality of the coefficients of the terms multiplying $\bar\beta
\bar\beta, \bar\beta$, and $1$. The $\bar\beta\bar\beta$ yield the three $\Phi_n^*,\Psi_n^*$
Poisson brackets and the terms with no $\beta$ yield the three $\Phi_n,\Psi_n$ Poisson brackets.
\end{proof}

\begin{remarks} 1. The no $\beta$ terms involve $z\Phi_n(z)$ and $w\Psi_n(w)$, but
$zw$ factors outs. The change of sign between \eqref{11.24} and \eqref{11.25} comes from
the minus sign in \eqref{4.1} vs.\ the plus sign in \eqref{4.2}.

\smallskip
2. It is an interesting exercise to use $\Phi_n^*(z) = z^n \, \ol{\Phi_n (1/\bar z)}$ to
go from \eqref{11.25} to \eqref{11.24}. Since Poisson brackets of real functions are
real, $\{\bar f, \bar g\}=\ol{\{f,g\}}$. We get a minus sign from $\ol{-i/2} =i/2$,
explaining the sign changes in the second and third terms. Because
$\ol{(\f{1/\bar z + 1/\bar w}{1/\bar z - 1/\bar w})} = -\f{z+w}{z-w}$,
the first term does not change.
\end{remarks}

\begin{theorem}\lb{T11.4} We have that
\begin{equation} \lb{11.26}
\{\Phi_n(z), \Phi_n^*(w)\} = i\,w \biggl( \f{\Phi_n(z) \Phi_n^*(w) - \Phi_n(w) \Phi_n^*(z)}
{z-w}\biggr)
\end{equation}
\end{theorem}

\begin{remarks} 1. The same formula holds if $\Phi_n,\Phi_n^*$ is replaced by $\Psi_n,
\Psi_n^*$ (since $\alpha_n \to -\alpha_n$ preserves $\{\cdot,\cdot\}$).

\smallskip
2. The proof provides another proof of $\{\Phi_n(z), \Phi_n(w)\}=\{\Phi_n^*(z),
\Phi_n^*(w)\}=0$.
\end{remarks}

\begin{proof} Define
\begin{equation} \lb{11.27}
\ti C_n(z) = z\Phi_n(z) \qquad \ti S_n(z) = \Phi_n^*(z)
\end{equation}
Then Szeg\H{o} recursion becomes
\begin{align}
\ti C_{n+1}(z) &= z\ti C_n(z) -\bar\alpha_n z\ti S_n(z) \lb{11.28} \\
\ti S_{n+1}(z) &= \ti S_n(z) -\alpha_n \ti C_n(z) \lb{11.29}
\end{align}
which has the same structure as \eqref{4.16}, \eqref{4.17} except that $\alpha_0,
\alpha_1,\dots$ is replaced by $\bar\alpha_0, \bar\alpha_1, \dots, \bar\alpha_n$.
Moreover, at $n=0$,
\[
\ti C_n (z) w \ti S_n(w) - \ti C_n(w) z\ti S_n(z) = zw-zw =0
\]
So we can start the induction. Thus, by induction as in Section~\ref{s5} $(-i$ becomes
the $+i$ because of the complex conjugate change which flips signs of $\{\alpha_n,
\bar\alpha_n\}$ to $\{\bar\alpha_n, \alpha_n\}$),
\[
\{\ti C_n(z), \ti S_n(w)\} = i \biggl( \f{\ti C_n(z) w\ti S_n(w) - \ti C_n(w) z\ti S_n(z)}
{z-w}\biggr)
\]
which is equivalent to \eqref{11.26}.
\end{proof}

Using the recursion relations for $\Phi_n$ in terms of $\Phi_{n-1}$, one obtains

\begin{theorem}\lb{T13.6} We have
\begin{align}
\{\Phi_{n-1}(z), \Phi_n(w)\} &= -i\bar\alpha_{n-1}w
\biggl( \f{\Phi_{n-1}(z) \Phi_{n-1}^*(w) - \Phi_{n-1}(w)\Phi_{n-1}^*(z)}
{z-w}\biggr) \lb{13.31} \\
\{\Phi_n(z), \Phi_{n-1}^*(w)\} &= -izw
\biggl( \f{\Phi_{n-1}(z) \Phi_{n-1}^*(w) - \Phi_{n-1}(w) \Phi_{n-1}^*(z)}
{z-w}\biggr) \lb{13.32}
\end{align}
\end{theorem}


\end{document}